\documentclass[sn-mathphys,Numbered]{sn-jnl}

\usepackage{graphicx}%
\usepackage{multirow}%
\usepackage{amsmath,amssymb,amsfonts}%
\usepackage{amsthm}%
\usepackage{mathrsfs}%
\usepackage[title]{appendix}%
\usepackage{xcolor}%
\usepackage{textcomp}%
\usepackage{manyfoot}%
\usepackage{booktabs}%
\usepackage{algorithm}%
\usepackage{algorithmicx}%
\usepackage{algpseudocode}%
\usepackage{listings}%

\theoremstyle{thmstyleone}%
\newtheorem{theorem}{Theorem}
\newtheorem{proposition}[theorem]{Proposition}%
\newtheorem{corollary}[theorem]{Corollary}

\theoremstyle{thmstyletwo}%
\newtheorem{example}{Example}%
\newtheorem{remark}{Remark}%

\theoremstyle{thmstylethree}%
\newtheorem{definition}{Definition}%

\begin{document}


 \title[Potential Vector Fields in $\mathbb R^3$]
  {Potential Vector Fields in $\mathbb R^3$ and $\alpha$-Meridional Mappings of the Second Kind 
$(\alpha \in \mathbb R)$}



\author*{\fnm{Dmitry} \sur{Bryukhov}} \email{bryukhov@mail.ru https://orcid.org/0000-0002-8977-3282}
\affil*{
 \orgname{Independent scholar},
\orgaddress{\street{Mira Avenue 19, apt. 225}, \city{Fryazino},
\postcode{141190}, \state{Moscow region}, \country{Russian
Federation}}}

\abstract{This paper extends author's approach developed in a recent  paper  
on analytic models of potential fields in inhomogeneous isotropic media. 
The properties of different analytic models in Cartesian and cylindrical coordinates in $\mathbb R^3$ are compared. 
The specifics of the Jacobian matrix $\mathbf{J}(\vec V)$ of potential meridional fields $\vec V$
 in cylindrically layered media, where $\phi( \rho) = \rho^{-\alpha}$ $(\alpha \in \mathbb R)$, 
lead to the concept of $\alpha$-meridional mappings of the first and second kind.
As a consequence, a new mathematical concept of $\alpha$-meridional functions of the first and second kind arises.
The topological and geometric properties of sets of degenerate points of the Jacobian matrix
 of various potential meridional fields $\vec V$ are demonstrated explicitly.
When $\alpha =1$, we deal with the special concept of radially holomorphic functions and the radially holomorphic potentials in $\mathbb R^3$. 
Remarkable properties of the radially holomorphic potentials represented by a linear superposition
of the radially holomorphic exponential function $e^{\breve{\beta} x}$ $(\breve{\beta} \in \mathbb R)$ 
 and function $e^{\breve{\beta} x}$ with coefficient $I$ are demonstrated explicitly. 
New properties of the radially holomorphic potentials represented 
by certain radially holomorphic fourth-order polynomials with reduced-quaternion-valued coefficients are established.
}




\maketitle

\section{Introduction}
\label{sec:intro}
The rich variety of three-dimensional analytic and numerical models of potential vector fields 
$\vec V = (V_0, V_1, V_2)$, where $V_l = V_l (x_0,x_1,x_2)$ $(l = 0,1,2)$, in mathematical physics and continuum mechanics 
(see, e.g., \cite{BornWolf:2003,Carslaw,KhmKravOv:2010,Reddy:2018,Br:Hefei2020})
 may be investigated by means of the following first-order system with a variable $C^1$-coefficient $\phi= \phi(x_0,x_1,x_2)>0$:
\begin{gather}
\begin{cases}
 \mathrm{div} \, (\phi \  \vec V) =0,  \\[1ex]
    \mathrm{curl}{\ \vec V} =0.
\end{cases}
\label{potential-system-3}
\end{gather}
The Euclidean space $\mathbb R^3=\{(x_0, x_1,x_2)\}$ in this setting involves the longitudinal variable $x_0$,
 the cylindrical radial variable $\rho = \sqrt{x_1^2+x_2^2}$, the polar angle $\varphi= \mathrm{arccot}\frac{x_0}{\rho}$ 
and the azimuthal angle $\ \theta = \arccos \frac{x_1}{\rho}$.

 The scalar potential $h = h(x_0,x_1,x_2)$  in simply connected open domains $\Lambda \subset \mathbb R^3$,
 where $\vec V = \mathrm{grad} \ h$,  allows us to reduce every
$C^1$-solution of the system~\eqref{potential-system-3} to a $C^2$-solution of the continuity equation
\begin{gather}
\mathrm{div} \, ( \phi \ \mathrm{grad}{\ h}) = 0.
\label{Liouville-3}
\end{gather}

In particular, the potential vector field $\vec V$ in the context of \emph{Mathematical problems in geometrical optics} 
may be interpreted as the electric field strength $\vec E$, and 
the coefficient $\phi$ as the dielectric permittivity $\varepsilon$ (see, e.g., \cite {BornWolf:2003,KhmKravOv:2010}), respectively.

 The scalar potential $h$ in the context of the theory of \emph{Conduction of heat} may be interpreted as the steady state temperature $T$, 
and the coefficient $\phi$ as the thermal conductivity $\kappa$  (see, e.g., \cite {Carslaw,Br:Hefei2020}).

The potential vector field $\vec V$, satisfying relations $\vec V = \frac {d{\vec x}}{dt} = \mathrm{grad} \ h$, 
in the context of \emph{Mathematical problems in continuum mechanics} may be interpreted as the potential velocity field  in the case of a steady flow,
and the scalar potential $h$ as the velocity potential (see, e.g., \cite{KochinKibelRoze:1964,Sedov:1994,Acheson,WhiteXue:2021}).

The geometric properties of the symmetric Jacobian matrix $\mathbf{J}(\vec V)$ in three dimensions, 
where $ \mathbf{J_{l m}}(\vec V) = \frac{\partial{V_l}}{\partial{x_m}}$ $(l, m = 0,1,2)$, 
are difficult to treat in contrast to properties of the symmetric  Jacobian matrix in two dimensions
 into the framework of the concept of \emph{Conformal mappings of the second kind}
(see, e.g., \cite{KochinKibelRoze:1964,LavSh:1987,Acheson,WhiteXue:2021}).

It should be noted here that the system~\eqref{potential-system-3} 
under the condition $\phi(\rho) = \rho^{-\alpha}$ $(\rho >0)$  in the expanded form is described as
\begin{gather}
\begin{cases}
   \mathrm{div}\ { \vec V} -
 \alpha \left( \frac{x_1}{\rho^2} V_1 + \frac{x_2}{\rho^2} V_2 \right) =0,  \\[1ex]
     \mathrm{curl}{\ \vec V} =0.
\end{cases}
\label{alpha-axial-hyperbolic-system-3}
\end{gather}

The  corresponding continuity equation~\eqref{Liouville-3} is written as
 \begin{gather}
(x_1^2+x_2^2)\Delta{h} - \alpha \left(  x_1\frac{\partial{h}}{\partial{x_1}} + x_2\frac{\partial{h}}{\partial{x_2}}\right)  =0.
  \label{eq-axial-hyperbolic-3-alpha}
  \end{gather}

General class of $C^1$-solutions of the system~\eqref{alpha-axial-hyperbolic-system-3}
 in the context of \emph{Quaternionic analysis in $\mathbb R^3$ and its non-Euclidean modifications}
(see, e.g., \cite{Leut:2000,LeZe:CMFT2004,Br:Hefei2020}) is equivalently represented  as general class of $C^1$-solutions of
a family of axially symmetric generalizations of the Cauchy-Riemann system in $\mathbb R^3$
\begin{gather}
\begin{cases}
(x_1^2+x_2^2) \left( \frac{\partial{u_0}}{\partial{x_0}}-
      \frac{\partial{u_1}}{\partial{x_1}}-\frac{\partial{u_2}}{\partial{x_2}} \right)
      + \alpha (x_1u_1+x_2u_2)=0, \\[1ex]
      \frac{\partial{u_0}}{\partial{x_1}}=-\frac{\partial{u_1}}{\partial{x_0}},
       \quad \frac{\partial{u_0}}{\partial{x_2}}=-\frac{\partial{u_2}}{\partial{x_0}}, \\[1ex]
      \frac{\partial{u_1}}{\partial{x_2}}=\ \ \frac{\partial{u_2}}{\partial{x_1}},
\end{cases}
\label{A_3^alpha-system}
\end{gather}
where $(u_0, u_1, u_2)=(V_0, -V_1, -V_2)$.

  Certain advantages of three-dimensional analytic models of potential vector fields $\vec V$
 in special cylindrically layered media, where $\phi( \rho) = \rho^{-\alpha}$ $(\alpha \in \mathbb R)$, were characterized in 2021 \cite{Br:Hefei2020}
using cylindrical coordinates.

Potential meridional fields are provided by the condition $ \frac{\partial{h}}{\partial{\theta}} = 0$
(see, e.g., \cite{KhmKravOv:2010,Br:Hefei2020}).
Potential transverse fields are provided by the condition $\frac{\partial{h}}{\partial{x_0}} = 0$, respectively. 

Rich applications of $\alpha$-meridional mappings of the second kind  in mathematical physics
 in the context of potential meridional fields in cylindrically layered media, where $\phi( \rho) = \rho^{-\alpha}$, $\alpha \ge 0$,  
were first  studied in 2021 \cite{Br:Hefei2020}. 

The main goal of this paper is to demonstrate explicitly new applications of $\alpha$-meridional mappings of the second kind 
 in fluid mechanics of layered media.

The paper is organized as follows.
In Section 2, basic concepts of \emph{Reduced quaternion-valued functions} are described in the first subsection.
Basic concepts of \emph{Potential vector fields in $\mathbb R^3$} are described in the second subsection.
Basic concepts of \emph{Autonomous systems and gradient systems} 
are described in the third subsection.
In Section 3,  three-dimensional analytic models of potential velocity fields $\vec V$ 
in special layered media are constructed. Recent advances in the theory of boundary value problems 
for elliptic equations with singular coefficients in $\mathbb R^3$ are noted.
In Section 4, analytic models of potential meridional velocity fields in cylindrically layered media 
with the mass density  $\phi( \rho) = \rho^{-\alpha}$, where $\alpha \ge 0$, are studied. 
Applied properties of $\alpha$-meridional mappings of the second kind are characterized
in the context of \emph{Stability theory of gradient systems} in  $\mathbb R^3=\{(x_0, x_1,x_2)\}$.
In Section 5, properties of $1$-meridional mappings of the second kind  are considered
 in the context of \emph{Generalized axially symmetric potential theory (GASPT)}.
The concept of the radially holomorphic potential in $\mathbb R^3$ 
allow us to extend the concept of the complex potential within potential meridional velocity fields 
in cylindrically layered media with the mass density $\phi( \rho) = \rho^{-1}$. 
In Section 6, we conclude the paper by describing future work in the context of 
\emph{Quaternionic analysis in $\mathbb R^4$ and its non-Euclidean modifications}.

\section{Preliminaries}
 \label{sec2}

 \subsection{Reduced Quaternion-Valued Functions: Basic Concepts}
 \label{subsec21}

The real algebra of quaternions $\mathbb H$ is a four dimensional
skew algebra over the real field generated by real unity $1$. Three
imaginary unities $i, j,$ and $k$ satisfy to multiplication rules
\begin{gather*}
i^2 = j^2 = k^2  = ijk = -1, \quad ij = -ji = k.
\end{gather*}
The independent quaternionic variable is defined as
$$x = x_0 + ix_1 + jx_2  + kx_3.$$
The quaternion conjugation of $x$ is defined by the following
automorphism:
$$ x \mapsto \overline{x} := x_0 - ix_1 - jx_2 - kx_3.$$

If $\rho = \sqrt {x_1^2+x_2^2+x_3^2} > 0$, then  $x= x_0 + I \rho$, 
where  $ I = \frac{i x_1+ j x_2+ k x_3 }{\rho}$, $ I^2=-1.$

The independent quaternionic variable may be interpreted as the vector \\ 
$\vec x = (x_0, x_1, x_2, x_3)$ in $\mathbb R^4$, where we deal with the Euclidean norm  
$$
\| x \|^2 :=  x \overline{x} = x_0^2 + x_1^2 + x_2^2 + x_3^2 := r^2.
 $$

If $x_3 > 0$, the independent quaternionic variable in cylindrical coordinates in $\mathbb{R}^4$ 
is described as
 $x = x_0 + \rho (i\cos{\theta} + j \sin{\theta}\cos{\psi} + k\sin{\theta}\sin{\psi}),$ where

$x_1 = \rho \cos{\theta}, \quad x_2 = \rho \sin{\theta}\cos{\psi}$,  
$ \quad x_3 = \rho \sin{\theta}\sin{\psi},$

 $ \varphi=  \arccos \frac{x_0}{r} \ (0 < \varphi < \pi)$, 
 $\quad \theta = \arccos \frac{x_1}{\rho} \ (0 \leq \theta \leq 2\pi),$

$\psi = \mathrm{arccot} \frac{x_2}{x_3} \ (0 < \psi < \pi).$

The dependent quaternionic variable is defined as
$$
u = u_0 + iu_1 + ju_2 +  ju_3 \sim (u_0, u_1, u_2, u_3).
$$
The quaternion conjugation of $u$ is defined by the following
automorphism:
$$
u \mapsto \overline{u} := u_0 - iu_1 - ju_2 - ku_3.
$$

If  $x_3 = 0$, then we deal with the independent reduced quaternionic variable 
 $x = x_0 + ix_1 + jx_2.$
 The independent reduced quaternionic variable 
 may be interpreted as the vector $\vec x = (x_0, x_1, x_2)$ in $\mathbb R^3$.

 If $\rho > 0$, the independent reduced quaternionic variable in cylindrical coordinates 
in $\mathbb{R}^3$ is described as $x = x_0 + \rho (i\cos{\theta} + j \sin{\theta})$, where 

$\varphi=  \arccos \frac{x_0}{r} = \mathrm{arccot}\frac{x_0}{\rho} \ (0 < \varphi < \pi),
\quad \theta = \arccos \frac{x_1}{\rho} \ (0 \leq \theta \leq 2\pi).$

 The dependent reduced quaternionic variable is defined as
 $$
 u = u_0 + iu_1 + ju_2 \sim (u_0, u_1, u_2).
 $$

\begin{definition}
Let $\Omega \subset \mathbb R^3$ be an open set. Every continuously differentiable mapping
$u= u_0 + iu_1 + ju_2: \Omega \rightarrow \mathbb{R}^3$ is called 
the reduced quaternion-valued $C^1$-function in $\Omega$.
\end{definition}

Analytic models of three-dimensional harmonic potential fields
 $\vec V = \vec V(x_0,x_1,x_2)$ satisfy the Riesz system in $\mathbb R^3$
\begin{gather*}
\begin{cases}
       \mathrm{div}\ { \vec V} =0,  \\[1ex]
   \mathrm{curl}{\ \vec V} =0.
\end{cases}
\end{gather*}

General class of analytic solutions of the Riesz system in $\mathbb R^3$
 in the context of \emph{Quaternionic analysis in $\mathbb R^3$} 
(see, e.g., \cite{Leut:2000,BraDel:2003,Del:2007,GM:2011}) 
is equivalently represented as  general class of analytic solutions of the system
\begin{gather*}
(R)
\begin{cases}
  \frac{\partial{u_0}}{\partial{x_0}}-
      \frac{\partial{u_1}}{\partial{x_1}}- \frac{\partial{u_2}}{\partial{x_2}} =0,  \\[1ex]
       \frac{\partial{u_0}}{\partial{x_1}}=-\frac{\partial{u_1}}{\partial{x_0}},
     \quad \frac{\partial{u_0}}{\partial{x_2}}=-\frac{\partial{u_2}}{\partial{x_0}}, \\[1ex]
      \frac{\partial{u_1}}{\partial{x_2}}=\ \ \frac{\partial{u_2}}{\partial{x_1}},
\end{cases}
\end{gather*}
where  $(u_0, u_1, u_2):=(V_0, -V_1, -V_2)$.
Analytic solutions of the system $(R)$ are referred to as the reduced quaternion-valued monogenic functions 
$u= u_0 + iu_1 + ju_2$  with harmonic components $u_l= u_l(x_0,x_1,x_2)$ $(l= 0,1,2)$.

 Unfortunately, general class of reduced quaternion-valued monogenic functions does not include 
class of the reduced quaternionic power functions, where $u= u_0 + iu_1 + ju_2 = (x_0 + ix_1 + jx_2)^n$, 
$n \in \mathbb{Z}$ (see, e.g., \cite{Leut:CV20,Leut:2000}).

 A multifaceted analytic extension of the concept of the power series with real and complex coefficients
 in the context of  \emph{Modified quaternionic analysis in $\mathbb R^3$}
has been developed by Leutwiler and Eriksson-Bique since 1992 (see, e.g., \cite{Leut:CV17,Leut:CV20,Leut:Rud96,ErLe:1998}). 
An important concept of radially holomorphic functions in the context of the theory of \emph{Holomorphic functions in $n$-dimensional space} 
was introduced by G\"{u}rlebeck, Habetha and Spr\"{o}ssig in 2008 \cite{GuHaSp:2008}.

\subsection{Potential Vector Fields in $\mathbb R^3$ and the Scalar Potentials: Basic Concepts}
 \label{subsec22}
Numerous mathematical problems of two-dimensional analytic models of potential fields $\vec V = \vec V(x,y)$
in homogeneous media have been studied by means of the complex potential.
 In accordance with the theory of holomorphic functions of a complex variable, 
where $f = f(z) = u + iv$, $z = x + iy$ \cite{LavSh:1987,Br:Hefei2020}, analytic models of 
potential velocity fields $\vec V$ in continuum mechanics are characterized by the principal invariants 
\begin{gather*}
I_{\mathbf{J}(\vec V)} \equiv \mathrm{tr} \mathbf{J}(\vec V) = 0, \quad
II_{\mathbf{J}(\vec V)} \equiv \det\mathbf{J}(\vec V) = - \mid f'(z) \mid^2 \leq 0.
\end{gather*}

General class of $C^1$-solutions of the system ~\eqref{potential-system-3}
may be equivalently represented as general class of $C^1$-solutions of the system 
\begin{gather}
\begin{cases}
\phi  \left( \frac{\partial{u_0}}{\partial{x_0}} - \frac{\partial{u_1}}{\partial{x_1}} -
\frac{\partial{u_2}}{\partial{x_2}}\right) + \left(\frac{\partial{\phi}}{\partial{x_0}}u_0 - \frac{\partial{\phi}}{\partial{x_1}}u_1 - \frac{\partial{\phi}}{\partial{x_2}}u_2\right) =0,\\[1ex]
\frac{\partial{u_0}}{\partial{x_1}}=-\frac{\partial{u_1}}{\partial{x_0}}, \quad
\frac{\partial{u_0}}{\partial{x_2}}=-\frac{\partial{u_2}}{\partial{x_0}}, \\[1ex]
\frac{\partial{u_1}}{\partial{x_2}}=\frac{\partial{u_2}}{\partial{x_1}},
\end{cases}
\label{Bryukhov-Kaehler-3}
\end{gather}
where $ (u_0, u_1, u_2)=(V_0, -V_1, -V_2)$ \cite{Br:Hefei2020}. 

 The system~\eqref{Bryukhov-Kaehler-3} may be characterized as generalized non-Euclidean modification
 of the system $(R)$ with respect to the conformal metric  \cite{Br:Hefei2020}:
 \begin{gather}
 ds^2 = \phi^2 (d{x_0}^2 + d{x_1}^2 + d{x_2}^2).
\label{Riemannian conformal metric}
\end{gather}
  The system~\eqref{A_3^alpha-system} under the condition $\alpha>0$ 
may be characterized as $\alpha$-axial-hyperbolic non-Euclidean modification of the system $(R)$ 
with respect to the conformal metric~\eqref{Riemannian conformal metric} defined outside the axis $x_0$ by formula \cite{Br:Hefei2020}:
\begin{gather*}
ds^2 = \frac{d{x_0}^2 + d{x_1}^2 + d{x_2}^2}{\rho^{2\alpha}}.
\end{gather*}
\begin{definition}
Every exact solution of eqn~\eqref{eq-axial-hyperbolic-3-alpha} under the condition $\alpha>0$ 
 in a simply connected open domain $\Lambda \subset \mathbb R^3$ $(\rho > 0)$
 is called $\alpha$-axial-hyperbolic harmonic potential in $\Lambda$.
\end{definition}
The continuity equation~\eqref{Liouville-3} in the expanded form is expressed as
\begin{gather}
\phi \Delta h
 + \frac{\partial{\phi}}{\partial{x_0}} \frac{\partial{h}}{\partial{x_0}} 
+ \frac{\partial{\phi}}{\partial{x_1}} \frac{\partial{h}}{\partial{x_1}} +
  \frac{\partial{\phi}}{\partial{x_2}}\frac{\partial{h}}{\partial{x_2}} =0.
\label{Liouville-eq-3-expanded}
\end{gather}

\begin{definition} 
 A set of points in a simply connected open domain $\Lambda$ defined implicitly as
\begin{gather*}
 {\mathfrak{S}}_C = \{(x_0,x_1,x_2) \in \Lambda : h(x_0,x_1,x_2) - C = 0 \},
\end{gather*}
where $C$ is an arbitrary constant, is said to be a level set of a potential function $h$.
\end{definition}

Equation of the form
\begin{gather}
h_C(x_0,x_1,x_2)  = h(x_0,x_1,x_2) - C = 0
 \label{level-C-3}
\end{gather}
allows us to establish the main properties of level sets of $h$ in $\Lambda$,
where $C$ is viewed as a parameter (see, e.g., \cite{GoVoJoGa:2009,Zhang:2017}). 
Real-valued function $h_C = h_C(x_0,x_1,x_2)$ vanishes everywhere in $\Lambda$, in contrast to scalar potential $h = h(x_0,x_1,x_2)$.
This is an important property of $h_C$ for applying the Implicit Function Theorem  in $\Lambda$ (see, e.g., \cite{Nikolsky:1977,ZachThoe:1986}).

\begin{definition} 
 A surface in $\Lambda$ satisfying equation~\eqref{level-C-3} is called a level surface of $h$ in $\Lambda$.
\end{definition} 

\begin{definition} (see, e.g., \cite[p.6]{LlibreMurza:2018})
A level surface of $h$ is said to be regular in $\Lambda$ if $\mathrm{grad} \ h \neq 0$ in $\Lambda$. 
\end{definition} 

The Implicit Function Theorem  provides sufficient conditions under which~\eqref{level-C-3} may be solved for $x_0$, $x_1$ or $x_2$.
In particular,~\eqref{level-C-3} may be solved for $x_0$ in terms of $x_1$, $x_2$ and the constant $C$
in a neighborhood of any point $\vec {\breve{x}} = (\breve{x}_0, \breve{x}_1, \breve{x}_2) \in \Lambda$
at which
 $\frac{\partial{h_C(\breve{x}_0, \breve{x}_1, \breve{x}_2)}}{\partial{x_0}} = $
 $\frac{\partial{h(\breve{x}_0, \breve{x}_1, \breve{x}_2)}}{\partial{x_0}} \neq 0$.
The partial derivatives of the corresponding function $x_0 = x_0(x_1,x_2, C)$ with respect to $x_1$ and $x_2$ 
may be obtained using implicit differentiation:
\begin{gather*}
\frac{\partial{x_0}}{\partial{x_1}} = - \frac{\partial{h}}{\partial{x_1}} \left( \frac{\partial{h}}{\partial{x_0}} \right)^{-1}, \quad
\frac{\partial{x_0}}{\partial{x_2}} = - \frac{\partial{h}}{\partial{x_2}} \left( \frac{\partial{h}}{\partial{x_0}} \right)^{-1}.
\end{gather*}

Using the total differential $dh$, eqn~\eqref{level-C-3} may be reformulated as an exact differential equation 
\begin{gather*}
 dh  = \frac{\partial{h}}{\partial{x_0}} d{x_0} + \frac{\partial{h}}{\partial{x_1}} d{x_1} 
+ \frac{\partial{h}}{\partial{x_2}} d{x_2} = 0.
\end{gather*}

\begin{definition} (see, e.g., \cite[p.61]{ArnoldGeom})
Let $\varsigma$ be a real-valued independent variable.
Assume that a first-order linear homogeneous partial differential equation 
\begin{gather}
  W_0 \frac{\partial{h}}{\partial{x_0}} + W_1 \frac{\partial{h}}{\partial{x_1}} + W_2 \frac{\partial{h}}{\partial{x_2}} = 0
\label{PDE}
\end{gather}
 is satisfied in $ \Lambda$ such that 
\begin{gather}
 \frac{dx_l}{d\varsigma} = W_l(x_0,x_1,x_2) \quad (l = 0,1,2),
\label{characteristics-3}
\end{gather}
 where $W_l(x_0,x_1,x_2) \in C^1(\Lambda)$.
Then autonomous ordinary differential equations~\eqref{characteristics-3} are called 
the characteristic equations of the partial differential equation~\eqref{PDE}, 
and the vector field $\vec W  = (W_0, W_1, W_2)$ is called the characteristic vector field of~\eqref{PDE}, respectively.
\end{definition}

 When $\vec W_l(\vec x) \neq 0$ $(l = 0,1,2)$ in $\Lambda$, the characteristic equations of~\eqref{PDE} may be written in \emph{symmetric form}: 
  \begin{gather*}
   \frac{dx_0}{W_0} = \frac{dx_1}{W_1} = \frac{dx_2}{W_2}.
  \end{gather*}

   According to \cite[pages 77-78]{BerGiac:2003}, the elementary expression $\vec W  = (W_0, W_1, W_2)$ of the vector field as well as its operator version 
 $ \mathcal{W} =  W_0 \frac{\partial{}}{\partial{x_0}} + W_1 \frac{\partial{}}{\partial{x_1}} + W_2 \frac{\partial{}}{\partial{x_2}}$
   is employed. 

 As noted by Nucci in the general case of $\mathbb R^n$ \cite[p.286]{Nucci:2005},
``The method of the Jacobi last multiplier provides a means to determine all the solutions of the partial differential equation
 \begin{gather*}
 \mathcal{A}f = \sum_{i= 1}^n  a_i(x_1, \ldots, x_n) \frac{\partial{f}}{\partial{x_i}}  = 0
\end{gather*}
or its equivalent associated Lagrange’s system
\begin{gather*}
 \frac{dx_1}{a_1} = \frac{dx_2}{a_2} = \ldots = \frac{dx_n}{a_n}."
\end{gather*}

Equation~\eqref{PDE} is geometrically characterized as the orthogonality condition 
for  the characteristic vector field $\vec W(\vec x) = \frac{d \vec x}{d\varsigma}$
and the gradient vector field $\vec V(\vec x) = \mathrm{grad} \ h(\vec x)$:
\begin{gather*}
  ( \vec W , \vec V )  = 0.
\end{gather*}

\begin{definition}(see, e.g., \cite[p.29, p.44]{ZachThoe:1986})
 Let $\vec W(\vec x)$ and $\vec V(\vec x)$ be nonvanishing vector fields in $\Lambda$.
 A $C^2$ scalar potential $h = h(\vec x)$ is called a first integral of the characteristic equations~\eqref{characteristics-3} of~\eqref{PDE} in $\Lambda$
if  $\vec W(\vec x)$ is orthogonal to $\mathrm{grad} \ h(\vec x)$ in $\Lambda$.
\end{definition}

\begin{definition}
A point $\vec x^* = (x_0^*,x_1^*,x_2^*) \in \Lambda$  is said to be a critical point of the scalar potential $h$ 
if $ \mathrm{grad} \ h(x_0^*,x_1^*,x_2^*) =0$.
The set of all critical points is called the critical set of $h$ in $\Lambda$, respectively.
 \end{definition}
\begin{remark}
As follows from three conditions
$\frac{\partial{h(x_0^*,x_1^*,x_2^*)}}{\partial{x_0}}  =0$,
$\frac{\partial{h(x_0^*,x_1^*,x_2^*)}}{\partial{x_1}}  =0$,
$\frac{\partial{h(x_0^*,x_1^*,x_2^*)}}{\partial{x_2}} =0$,
eqn~\eqref{Liouville-eq-3-expanded} takes a simplified form $ \Delta h =0$ within the critical set of $h$.
  \end{remark}
\begin{definition}
Let $\mathbf{H}(h)$ be the Hessian matrix of $h$.
A critical point $\vec x^* = (x_0^*,x_1^*,x_2^*) \in \Lambda$ of the scalar potential $h = h(x_0, x_1, x_2)$ 
is said to be a non-Morse (or \emph{degenerate}, or \emph{nonisolated}) critical point if $\det\mathbf{H}(h(x_0^{*},x_1^{*},x_2^{*})) =0$.
 Otherwise, it is called a Morse (or \emph{nondegenerate}, or \emph{isolated}) critical point of $h$.
 \end{definition}

 \begin{remark}
 It is well known that only non-Morse critical points exist 
within the framework of the theory of \emph{Holomorphic functions of a complex variable} (see e.g., \cite{LavSh:1987}).
   \end{remark}

The characteristic equation of the Jacobian matrix of arbitrary potential $C^1$-vector field $\vec V$
in the general setting
\begin{gather}
\begin{pmatrix}
  \frac{\partial{V_0}}{\partial{x_0}} & \frac{\partial{V_0}}{\partial{x_1}} & \frac{\partial{V_0}}{\partial{x_2}} \\[1ex]
 \frac{\partial{V_1}}{\partial{x_0}}  & \frac{\partial{V_1}}{\partial{x_1}}  &  \frac{\partial{V_1}}{\partial{x_2}}  \\[1ex]
 \frac{\partial{V_2}}{\partial{x_0}}  & \frac{\partial{V_2}}{\partial{x_1}}  &  \frac{\partial{V_2}}{\partial{x_2}}
 \end{pmatrix} =
\begin{pmatrix}
 \ \ \frac{\partial{u_0}}{\partial{x_0}} &  \ \ \frac{\partial{u_0}}{\partial{x_1}} & \ \ \frac{\partial{u_0}}{\partial{x_2}}  \\[1ex]
 -\frac{\partial{u_1}}{\partial{x_0}}  & -\frac{\partial{u_1}}{\partial{x_1}}  &  -\frac{\partial{u_1}}{\partial{x_2}}  \\[1ex]
 -\frac{\partial{u_2}}{\partial{x_0}}  & -\frac{\partial{u_2}}{\partial{x_1}}  &  -\frac{\partial{u_2}}{\partial{x_2}}
 \end{pmatrix}
 \label{Hessian-matrix-3}
\end{gather}
  is expressed as (see e.g., \cite{LaiRubKr:2010,Br:Hefei2020}) 
\begin{gather}
\lambda^3 - I_{\mathbf{J}(\vec V)} \lambda^2 + II_{\mathbf{J}(\vec
V)} \lambda - III_{\mathbf{J}(\vec V)}  = 0.
  \label{characteristic lambda-3}
\end{gather}

The principal scalar invariants $I_{\mathbf{J}(\vec V)}$, $II_{\mathbf{J}(\vec V)}$, $III_{\mathbf{J}(\vec V)}$ 
 are given by the formulas
\begin{gather}
\begin{cases}
I_{{\mathbf{J}(\vec V)}} \equiv \mathrm{tr} \mathbf{J}(\vec V) \equiv  \mathrm{div} \, \vec V = \lambda_0 +
\lambda_1 + \lambda_2= J_{00} + J_{11} + J_{22},  \\[1ex]
II_{{\mathbf{J}(\vec V)}} = \lambda_0 \lambda_1 + \lambda_0 \lambda_2 + \lambda_1 \lambda_2 = \\[1ex]
 J_{00}J_{11} + J_{00}J_{22} + J_{11}J_{22} - (J_{01})^2 - (J_{02})^2 - (J_{12})^2, \\[1ex]
III_{{\mathbf{J}(\vec V)}} \equiv \det\mathbf{J}(\vec V) =  \lambda_0 \lambda_1 \lambda_2 = \\[1ex]
  J_{00}J_{11}J_{22} + 2J_{01}J_{02}J_{12} - J_{00}(J_{12})^2 - J_{11}(J_{02})^2 - J_{22}(J_{01})^2,
\end{cases}
 \label{principal invariants}
\end{gather}
where real roots $\lambda_0$, $\lambda_1$, $\lambda_2$ of eqn~\eqref{characteristic lambda-3} 
are the eigenvalues of~\eqref{Hessian-matrix-3}.

 The principal scalar invariants~\eqref{principal invariants} in $\mathbb R^3$ play key roles within 
analytic models of potential fields in mathematical physics and continuum mechanics
 (see, e.g., \cite{LaiRubKr:2010,Br:Hefei2020}).
 The third principal invariant may have a variable sign in simply connected open domains $\Lambda \subset \mathbb R^3$
 in contrast to the second principal invariant within the framework of the concept of \emph{Conformal mappings of the second kind}
 (see e.g., \cite{KochinKibelRoze:1964,LavSh:1987,Acheson,WhiteXue:2021,Br:Hefei2020}). 

\begin{remark}
According to \cite[p.20]{PosStew}, a square matrix $\mathbf{A}$ is non-singular if and only if $\det\mathbf{A} \neq 0$.
According to \cite[Basic Concepts]{ArnoldGusVar:2012}, the theory of singularities of smooth maps 
is connected with the vanishing of certain derivatives and Jacobians.
  \end{remark}

To avoid conflicting terminology in the context of the modern theory of \emph{Elliptic equations with singular coefficients} 
(see, e.g., \cite{KarimETNieto:2011,NietoKarimET:2012,UrinovKarimovKT:2019,UriKar:2023}),
the following definition may be introduced (see, e.g., \cite{Br:Hefei2020}).

  \begin{definition} 
   A point $\vec x = (x_0,x_1,x_2) \in \Lambda$ is said to be a degenerate point  of the Jacobian matrix $\mathbf{J}(\vec V)$
 in $\Lambda$ if $\det\mathbf{J}(\vec V(\vec x)) =0$.
Otherwise, it is called a nondegenerate point of $\mathbf{J}(\vec V)$ in $\Lambda$.
   \end{definition}

The Jacobian matrix $\mathbf{J}(\vec V)$ of arbitrary potential $C^1$-vector field $\vec V$ 
coincides with the Hessian matrix $\mathbf{H}(h)$ of the corresponding scalar potential $h$.
 Sets of degenerate points of $\mathbf{J}(\vec V)$ of potential fields $\vec V=  \mathrm{grad} \ h$ in $\Lambda \subset \mathbb R^3$
   cover sets of non-Morse critical points of potential functions $h$ in $\Lambda \subset \mathbb R^3$.

 \begin{remark}
 The topological and geometric properties of sets of degenerate points of the Jacobian matrix of harmonic vector fields in $\mathbb R^3$ are of particular interest. 
 Following the concept expounded by H. Lewy \cite{Lewy:1968}, Yanushauskas established in 1980 
 some surprising properties of zero eigenvalues of the Hessian matrix of harmonic functions in $\mathbb{R}^n$  $(n \geq 3)$ \cite{Yanush:1980}.
 The Hessian of harmonic functions assumes both positive and negative values in any sufficiently small neighborhood of critical points.
  \end{remark}

 \subsection
{Vector Fields in the Phase Space, Autonomous Systems and Gradient Systems: Basic Concepts}
 \label{subsec23}

To characterize the main properties of axisymmetric gradient systems in $\mathbb R^4=\{(x_0, x_1, x_2, x_3)\}$, 
it is necessary to delve into the theory of \emph{Autonomous systems of ordinary differential equations} 
in $\mathbb R^n=\{(x_1, \ldots, x_n)\}$
(see, e.g., \cite{NemStep:1960,Goriely:2001,Perko:2001,Wiggins:2003,BerGiac:2003,Zhang:2017,Strogatz:2018,Kuznetsov:2023}).
The Euclidean space $\mathbb R^n$ is referred to as \emph{the phase space} 
(see, e.g., \cite[p.1]{Wiggins:2003}), or \emph{the state space} (see, e.g., \cite[p.5]{Kuznetsov:2023}).
 
\begin{definition} (see, e.g., \cite[p.1]{Zhang:2017})
Let $\vec V = (V_1, \ldots, V_n)$ be a nonlinear vector field in an open set $\Omega \subset \mathbb R^n$.
An autonomous system of first-order ordinary differential equations 
\begin{gather}
  \frac{d \vec x}{dt} = \vec V(\vec x)
 \label{auton-n}
\end{gather}
 is said to be a smooth nonlinear dynamical system if $\vec V(\vec x) \in C^1(\Omega)$.
 \end{definition}

 \begin{definition} (see, e.g., \cite[p.440]{Llibre:2004})
 The real-valued independent variable $t$ is called the time of system~\eqref{auton-n}.
  \end{definition}

\begin{remark}
Nonlinear autonomous systems of first-order ordinary differential equations 
are often treated as \emph{Nonlinear differential systems} in the modern literature 
(see, e.g., \cite{Llibre:2004,DumLliArt:2006,LlibreZhang:2012,Zhang:2016,Zhang:2017}).
 \end{remark}

\begin{definition} 
A point $\vec x^{**} = (x_1^{**}, \ldots, x_n^{**}) \in \Omega$ is said to be an equilibrium point
of a smooth system~\eqref{auton-n} if $\vec V(\vec x^{**}) = 0$.
Otherwise, it is called a regular point of~\eqref{auton-n}.
The set of all equilibrium points in $\Omega$ is called the set of equilibria of~\eqref{auton-n} in $\Omega$, respectively.
 \end{definition}

\begin{definition} (see, e.g., \cite[p.78]{BerGiac:2003})
A $C^1$ real-valued function $M = M(\vec x) >0$ defined in $\Omega \subset \mathbb R^n$
is said to be a Jacobi multiplier for~\eqref{auton-n} in $\Omega$ 
when $M$ solves the following linear first-order partial differential equation:
\begin{gather*}
 \sum_{i= 1}^n  \frac{\partial{(M V_i)}}{\partial{x_i}}  = 0.
\end{gather*}
\end{definition}

 As noted by Berrone and Giacomini \cite[p.81]{BerGiac:2003}, 
``Jacobi multipliers have been also named ``last multipliers" in the classical literature"
(for more details on this see \cite{BerGiac:2003,Nucci:2005,Zhang:2017}).

  From a geometrical point of view, an arbitrary solution $\vec x = \vec x(t)$ of~\eqref{auton-n}
is treated as \emph{the trajectory}, or \emph{phase curve} in the phase space $\mathbb R^n$  (see, e.g., \cite[p.2]{Wiggins:2003}). 

\begin{definition} (see, e.g., \cite[p.7]{Zhang:2017})
A point $\vec x^{**} = (x_1^{**}, \ldots, x_n^{**}) \in \Omega$ is said to be an equilibrium point
of~\eqref{auton-n} if $\vec V(\vec x)  = 0$ at $\vec x^{**}$.
Otherwise, it is called a regular point of~\eqref{auton-n}.
Set of all equilibrium points in $\Omega$ is called the set of equilibria of~\eqref{auton-n} in $\Omega$, respectively.
 \end{definition}

\begin{remark}
The sets of equilibria of~\eqref{auton-n} in $\Omega$ coincide with zero sets of the vector fields $\vec V$ in $\Omega$.
 \end{remark}

\begin{definition} (see, e.g., \cite[p.102]{Perko:2001})
An equilibrium point $\vec x^{**} = (x_1^{**}, \ldots, x_n^{**}) \in \Omega$ of~\eqref{auton-n}
 is said to be a hyperbolic point if none of the eigenvalues 
of the Jacobian matrix $\mathbf{J}(\vec V(\vec x))$ at $\vec x^{**}$ have zero real part.
Otherwise, it is called a nonhyperbolic point of~\eqref{auton-n}.
   \end{definition}

 \begin{definition}
  An equilibrium point $\vec x^{**} \in \Omega$ of~\eqref{auton-n} is said to be a degenerate if $\det\mathbf{J}(\vec V(\vec x^{**})) =0$.
 Otherwise, it is called a nondegenerate equilibrium point of~\eqref{auton-n}.
   \end{definition}

Equilibrium points of~\eqref{auton-n} in the context of \emph{Qualitative theory of differential equations}, 
the theory of \emph{Integrability}, \emph{Stability theory} and \emph{Bifurcation theory} are very often viewed as  \emph{singular points},
or \emph{zeros}, or \emph{critical points}, or \emph{fixed  points}, or \emph{stationary points}
(see, e.g., \cite{NemStep:1960,Chetayev:1990,Perko:2001,Wiggins:2003,DumLliArt:2006,Strogatz:2018,Goriely:2001,Zhang:2017,Kuznetsov:2023}).

As noted by Strogatz \cite[p.47]{Strogatz:2018}, ``Bifurcation theory is rife with conflicting terminology.
     The subject really hasn't settled down yet, and different people use different words for the same thing."

\begin{definition} 
Let $h = h(\vec x)$ be a $C^2$ single-valued scalar potential in a simply connected open domain $\Lambda \subset \mathbb R^n$,
 such that $\vec V(\vec x) = \mathrm{grad} \ h(\vec x)$. 
A smooth dynamical system of the form 
\begin{gather}
  \frac{d \vec x}{dt} = \vec V(\vec x) = \mathrm{grad} \ h(\vec x)
 \label{grad-n}
\end{gather}
 is said to be a smooth gradient system with scalar potential $h$ in $\Lambda$.
 \end{definition}

\begin{remark}
 When we deal with gradient systems~\eqref{grad-n}, all eigenvalues of the Jacobian matrix $\mathbf{J}(\vec V(\vec x))$ 
at an arbitrary equilibrium point $\vec x^{**}$ have zero imaginary parts, and
 sets of nonhyperbolic points coincide with sets of degenerate equilibrium points.
 \end{remark}

 Consider dynamical systems in the phase space $\mathbb R^n=\{(x_1, \ldots, x_n)\}$ in a broader context, 
where $C^1$ vector fields $\vec V = (V_1, \ldots, V_n)$ depend on several parameters $\mu_1 \ldots, \mu_{\check{l}}$
(see, e.g., \cite{ArnAfrIlyashShil:1994,Kuznetsov:2023}):
\begin{gather}
  \frac{d \vec x}{dt} = \vec V(\vec x; \mu_1, \ldots, \mu_{\check{l}}).
 \label{auton-n-mu}
\end{gather}
Dynamical systems of the form~\eqref{auton-n-mu} are known as dynamical systems depending on several parameters 
(see, e.g., \cite{Perko:2001,HirschSmaleDev:2013,Kuznetsov:2023}).
The Euclidean space $\mathbb R^{\check{l}} =\{(\mu_1, \ldots, \mu_{\check{l}})\}$
is known as  the space of parameters (see, e.g.,  \cite{ArnAfrIlyashShil:1994}).

According to \cite[pages 44-45]{Chetayev:1990}, the potential energy of material systems in mechanics often depends on several parameters.
When the Hessian of the potential energy vanishes for certain fixed parameter values, bifurcations of equilibria occur.

 Similarly, bifurcations of equilibria of a smooth gradient system
\begin{gather*}
  \frac{d \vec x}{dt} = \vec V(\vec x; \mu_1, \ldots, \mu_{\check{l}}) = -\mathrm{grad} \ \hslash(\vec x; \mu_1, \ldots, \mu_{\check{l}})
\end{gather*}
with potential function $\hslash = \hslash(\vec x; \mu_1, \ldots, \mu_{\check{l}})$
 for certain fixed parameter values $(\mu_1^{**}, \ldots, \mu_{\check{l}}^{**})$ occur at $x^{**}$ in $\Lambda$, 
when $\det\mathbf{H}(\hslash(\vec x^{**}; \mu_1^{**}, \ldots, \mu_{\check{l}}^{**})) =0$.

The advantages of studying \emph{Polynomial differential systems} have attracted increased interest from researchers
in the context of the theory of \emph{Integrability} 
(see, e.g., \cite{Singer:1992,Goriely:2001,LlibreZhang:2002,Llibre:2004,DumLliArt:2006,LlibreZhang:2012,Zhang:2017,LlibreMurza:2018}).
 As noted by Zhang \cite[p.306]{Zhang:2011}, ``the Darboux theory of integrability was developed 
for polynomial vector fields in $\mathbb C^n$ and $\mathbb R^n$".

According to \cite[pages 21-22]{Zhang:2017},
especially interesting cases of the characteristic equations~\eqref{characteristics-3} of~\eqref{PDE}
 may be investigated in the context of the theory of \emph{Integrability of polynomial differential systems in $\mathbb R^3$}. 

 Some remarkable properties of polynomial systems in $\mathbb R^4$ 
represented by the so-called one-dimensional quaternion homogeneous polynomial differential equation  
 \begin{gather}
  \frac{dq}{dt} = \check{a} q^{\check{k}}\overline{q}^{\check{n}},
 \label{a-overline-monomial-k,n}
\end{gather} 
where $\check{a} \in \mathbb H$,  $\check{k}, \check{n} \in \mathbb N \bigcup \{0\}$, 
$q = q_0 + q_1i + q_2j + q_3k$ and $\overline{q}$ is the quaternion conjugation of $q$,
 were considered by Gasull, Llibre and Zhang in 2009 \cite{GasLliZh:2009}). 
According to \cite{GasLliZh:2009}, the right-hand side of~\eqref{a-overline-monomial-k,n} is an unique monomial.

 When $\check{n}= 0$, the quaternion differential equation~\eqref{a-overline-monomial-k,n} is written as
\begin{gather}
  \frac{dq}{dt} = \check{a} q^{\check{k}}.
\label{monomial-k}
\end{gather} 
Certain important cases of eqn~\eqref{monomial-k}, where $\check{a} \in \mathbb H$, were studied.

 When $\check{k}= 0$, eqn~\eqref{a-overline-monomial-k,n} is written as 
\begin{gather}
  \frac{dq}{dt} = \check{a} \overline{q}^{\check{n}}.
\label{overline-monomial-n}
\end{gather} 
Certain important cases of eqn~\eqref{overline-monomial-n}, where $\check{a} \in \mathbb H$, were highlighted.

Several new kinds of polynomial autonomous systems in $\mathbb R^4$ 
represented by polynomial differential equations over the quaternions  
\begin{gather}
  \frac{dx}{dt} = P(x),
\label{WaZh-polynomial}
\end{gather} 
where $x = x_0 + x_1i + x_2j + x_3k$ and $P(x)$ is a quaternionic polynomial with complex coefficients, 
were studied by Zhang in 2011 \cite{Zhang:2011} and by Walcher and Zhang in 2021 \cite{WalZhang:2021}. 
Qualitative properties of equilibrium (or \emph{stationary}) points of polynomial autonomous systems represented by~\eqref{WaZh-polynomial} 
raise new issues for consideration in the context of \emph{Stability theory}.

 A smooth gradient system with scalar potential $h$ in a simply connected open domain 
$\Lambda \subset \mathbb R^3=\{(x_0, x_1,x_2)\}$ may be described as  
(see, e.g., \cite{Wiggins:2003,HirschSmaleDev:2013,Strogatz:2018,BrRhod:2013,BrRhod:2014})
\begin{gather}
 \frac {d{\vec x}}{dt} = \vec V = \mathrm{grad} \ h(\vec x), \quad t \in \mathbb R.
\label{grad-system-3}
\end{gather}

 The first applications of gradient systems~\eqref{grad-system-3} were provided in 2013-2014 \cite{BrRhod:2013,BrRhod:2014}.
Potential (often referred to as \emph{irrotational}  \cite{LaiRubKr:2010,BrKos:2012,BrRhod:2013})
meridional velocity fields $\vec V$ in special inhomogeneous isotropic media with the mass density $\phi = \rho^{-1}$
 were represented by the following reduced-quaternion-valued ordinary differential equation:
\begin{gather*}
  \frac {dx}{dt} =  V_0 + i V_1 + j V_2 = \overline{F}(x),
  \end{gather*}
where $x= x_0 + ix_1 + jx_2$, $\overline{F}(x) = u_0 - i u_1 - j u_2$ and
 $F(x) = \frac{\partial{h}}{\partial{x_0}} - i \frac{\partial{h}}{\partial{x_1}} - j\frac{\partial{h}}{\partial{x_2}}$.

\section
{Analytic Models of Potential Velocity Fields in Some Special Layered Media}
 \label{sec3}

Inhomogeneous isotropic media, whose properties are constant throughout every plane perpendicular to a fixed direction, 
are referred in mathematical physics and continuum mechanics to as layered media (see, e.g., \cite {BornWolf:2003,Brekh:1980,Br:Hefei2020}).

Let us pay attention to remarkable properties of analytic models of potential velocity fields $\vec V$ 
in special planarly layered media, where $\phi = \phi_1(x_1)\phi_2(x_2)$, $\phi_1(x_1) >0$, $\phi_2(x_2) >0$:
\begin{gather}
\begin{cases}
  \mathrm{div} \, ( \phi_1(x_1)\phi_2(x_2) \vec V ) = 0, \\[1ex]
  \mathrm{curl}{\ \vec V} = 0.
\end{cases}
 \label{bi-potential-system-3}
\end{gather}

General class of $C^1$-solutions of the system~\eqref{bi-potential-system-3} is equivalently represented as general class of $C^1$-solutions of the system
\begin{gather}
\begin{cases}
      \phi_1(x_1)\phi_2(x_2) \left(\frac{\partial{u_0}}{\partial{x_0}}-
      \frac{\partial{u_1}}{\partial{x_1}}-  \frac{\partial{u_2}}{\partial{x_2}}\right)
  - \left( \frac{d{{\phi}_1}}{d{x_1}}u_1 + \frac{d{{\phi}_2}}{d{x_2}}u_2 \right)  = 0, \\[1ex]
      \frac{\partial{u_0}}{\partial{x_1}}=-\frac{\partial{u_1}}{\partial{x_0}}, \quad 
      \frac{\partial{u_0}}{\partial{x_2}}=-\frac{\partial{u_2}}{\partial{x_0}}, \\[1ex]
      \frac{\partial{u_1}}{\partial{x_2}}=\frac{\partial{u_2}}{\partial{x_1}},
\end{cases}
\label{Bryukhov-3-hyperbolic-3}
\end{gather}
where  $(V_0,V_1,V_2) = (u_0, -u_1, -u_2)$.

Eqn~\eqref{Liouville-eq-3-expanded} is written as
\begin{gather}
 \phi_1(x_1)\phi_2(x_2) \left( \frac{{\partial}^2{h}}{{\partial{x_0}}^2} + \frac{{\partial}^2{h}}{{\partial{x_1}}^2} + \frac{{\partial}^2{h}}{{\partial{x_2}}^2} \right)
 +   \frac{d{{\phi}_1}}{d{x_1}} \frac{\partial{h}}{\partial{x_1}} +
  \frac{d{{\phi}_2}}{d{x_2}} \frac{\partial{h}}{\partial{x_2}} =0.
\label{alpha_1,2-biplanar}
\end{gather}
  Suppose that  $\phi_1(x_1) =  x_1^{-\alpha_1}$, $\phi_2(x_2) =  x_2^{-\alpha_2}$ $(\alpha_1, \alpha_2 \in \mathbb{R})$.
Eqn~\eqref{alpha_1,2-biplanar} is reduced to the following elliptic equation with two singular coefficients: 
\begin{gather}
\Delta{h} - \frac{\alpha_1}{x_1}\frac{\partial{h}}{\partial{x_1}} - \frac{\alpha_2}{x_2}\frac{\partial{h}}{\partial{x_2}} =0.
\label{alpha_1,2-bihyperbolic-3}
\end{gather}

 The system~\eqref{bi-potential-system-3} is expressed as
\begin{gather*}
\begin{cases}
  \mathrm{div} \, ( x_1^{-\alpha_1} x_2^{-\alpha_2} \vec V ) = 0, \\[1ex]
  \mathrm{curl}{\ \vec V} = 0,
\end{cases}
\end{gather*}
and the system~\eqref{Bryukhov-3-hyperbolic-3}  is simplified:
\begin{gather*}
\begin{cases}
      (\frac{\partial{u_0}}{\partial{x_0}}-
      \frac{\partial{u_1}}{\partial{x_1}}-\frac{\partial{u_2}}{\partial{x_2}}) + \frac{\alpha_1}{x_1} u_1 + \frac{\alpha_2}{x_2} u_2 = 0, \\[1ex]
      \frac{\partial{u_0}}{\partial{x_1}}=-\frac{\partial{u_1}}{\partial{x_0}},  \quad 
     \frac{\partial{u_0}}{\partial{x_2}}=-\frac{\partial{u_2}}{\partial{x_0}}, \\[1ex]
      \frac{\partial{u_1}}{\partial{x_2}}=\ \ \frac{\partial{u_2}}{\partial{x_1}}.
\end{cases}
\end{gather*}
This system under conditions of $\alpha_1>0$, $\alpha_2>0$ may be characterized 
as $(\alpha_1, \alpha_2)$-bihyperbolic non-Euclidean modification of the system $(R)$
 with respect to the conformal metric~\eqref{Riemannian conformal metric} defined on a quarter-space $\{x_1 > 0, x_2 > 0\}$ by formula:
\begin{gather*}
ds^2 = \frac{d{x_0}^2 + d{x_1}^2 + d{x_2}^2}{ x_1^{2\alpha_1} x_2^{2\alpha_2}}.
\end{gather*}

 \begin{definition}
 Every exact solution of eqn~\eqref{alpha_1,2-bihyperbolic-3} under the conditions $\alpha_1>0$, $\alpha_2> 0$
 in a simply connected open domain $\Lambda \subset \mathbb R^3$ $(x_1 > 0, x_2 > 0)$ 
 is called $(\alpha_1, \alpha_2)$-bihyperbolic harmonic potential in $\Lambda$.
 \end{definition}

 The specifics of $(\alpha_1, \alpha_2)$-bihyperbolic harmonic potentials may be demonstrated explicitly using separation of variables. 
\begin{theorem}
A special class of three-dimensional solutions of eqn~\eqref{alpha_1,2-bihyperbolic-3} 
may be obtained using the Bessel functions of the first and second kind, trigonometric functions, power functions 
and the modified Bessel functions of the first and second kind for different values of the separation constants $\breve{\lambda}$ and $\breve{\mu}$:
\begin{align*}
& h(x_0, x_1, x_2) = {x_1}^\frac{\alpha_1+1}{2} \left[ c_{\breve{\lambda}}^1 J_{\frac{\alpha_1+1}{2}}(\breve{\lambda}x_1)
+ c_{\breve{\lambda}}^2 Y_{\frac{\alpha_1+1}{2}}(\breve{\lambda}x_1)  \right] \times \\
& \sum_{\breve{\mu}= -\infty}^\infty \left( b^1_{\breve{\mu}} \cos{\breve{\mu} x_0} +  b^2_{\breve{\mu}} \sin{\breve{\mu} x_0} \right)
{x_2}^\frac{\alpha_2+1}{2} \left[ a^1_{\breve{\lambda}, \breve{\mu}} I_{\frac{\alpha_2+1}{2}}( \breve{\nu}x_2) +
a^2_{\breve{\lambda}, \breve{\mu}} K_{\frac{\alpha_2+1}{2}}( \breve{\nu}x_2) \right],
\end{align*}
where $\ \breve{\nu} = \sqrt{ \breve{\lambda}^2 + \breve{\mu}^2}$;
$\ c^1_{\breve{\lambda}},  c^2_{\breve{\lambda}},
b^1_{\breve{\mu}}, b^2_{\breve{\mu}}, a^1_{\breve{\lambda}, \breve{\mu}}, a^2_{\breve{\lambda}, \breve{\mu}}
 = const \in \mathbb R $.
\end{theorem}
 \begin{proof}
 Consider a special class of exact solutions of eqn~\eqref{alpha_1,2-bihyperbolic-3} under the assumption that  $h(x_0, x_1, x_2) =$ $p(x_0,  x_2) \varpi(x_1)$:
$$
 \varpi \left( \frac{\partial{^2}{p}}{\partial{x_0}^2} +
         \frac{\partial {^2}{p}}{\partial{ x_2}^2} \right) -
\frac{\varpi \alpha_2}{x_2} \frac{\partial{p}}{\partial{ x_2}}  +
        p \frac{d{^2}{\varpi}}{d{x_1}^2} - \frac{ \alpha_1}{x_1} p \frac{d{\varpi}}{d{x_1}}
 = 0.
$$
 Relations
\begin{align*}
  - p \frac{d{^2}{\varpi}}{d{x_1}^2} + \frac{ \alpha_1}{x_1} p \frac{d{\varpi}}{d{x_1}} =
  \varpi \left( \frac{\partial{^2}{p}}{\partial{x_0}^2} +
         \frac{\partial {^2}{p}}{\partial{x_2}^2} \right)
 - \frac{\varpi \alpha_2}{x_2} \frac{\partial{p}}{\partial{ x_2}}  =
     \breve{\lambda}^2 p\varpi  \quad  ( \breve{\lambda} = const \in \mathbb R )
\end{align*}
lead to the following system of equations:
\begin{gather}
\begin{cases}
   \frac{d{^2}{\varpi}}{d{x_1}^2} - \frac{\alpha_1}{x_1} \frac{d{\varpi}}{d{x_1}} + \breve{\lambda}^2 \varpi = 0, \\
   \frac{\partial{^2}{p}}{\partial{x_0}^2} +   \frac{\partial {^2}{p}}{\partial{x_2}^2}
 - \frac{\alpha_2}{x_2} \frac{\partial{p}}{\partial{x_2}}
  - \breve{\lambda}^2 p = 0.
 \end{cases}
  \label{Laplace-Beltrami equation, bi-sep-3}
\end{gather}
The first equation of the system~\eqref{Laplace-Beltrami equation, bi-sep-3} 
as a linear second-order ordinary differential equation containing power functions 
may be solved using linear independent solutions 
(see, e.g., \cite[Chapter 14, p.526 item 63]{PolZait:Ordin-2018}):
$$ \varpi_{ \breve{\lambda}}(x_1)=
{x_1}^\frac{\alpha_1+1}{2} \left[ c_{\breve{\lambda}}^1 J_{\frac{\alpha_1+1}{2}}(\breve{\lambda}x_1)
+ c_{\breve{\lambda}}^2 Y_{\frac{\alpha_1+1}{2}}(\breve{\lambda}x_1)  \right]; \quad
  c_{\breve{\lambda}}^1, c_{\breve{\lambda}}^2= const \in \mathbb{R}, $$
 where  $J_{ \breve{\nu}}(\breve{\xi})$ and $Y_{ \breve{\nu}}(\breve{\xi})$ are the Bessel functions of the first and second kind of 
real order ${\frac{\alpha_1 + 1}{2}}$ and real argument $\breve{\lambda}x_1$ (see, e.g., \cite{Watson:1944,Koren:2002}).

The second equation of the system~\eqref{Laplace-Beltrami equation, bi-sep-3}
may be solved using separation of variables $p(x_0,  x_2) = \Xi(x_0) \Upsilon(x_2)$:
$$
\frac{1}{\Xi}  \frac{d{^2}{\Xi}}{d{x_0}^2} +   \frac{1}{ \Upsilon} \frac{d{^2}{ \Upsilon}}{d{x_2}^2}
- \frac{\alpha_2} { \Upsilon x_2} \frac{d{ \Upsilon}}{d{x_2}} -  \breve{\lambda}^2= 0.
$$
Relations
\begin{align*}
- \frac{1}{\Xi} \frac{d{^2}{\Xi}}{d{x_0}^2} =  \frac{1}{ \Upsilon} \frac{d{^2}{ \Upsilon}}{d{x_2}^2}
 - \frac{\alpha_2} { \Upsilon x_2} \frac{d{ \Upsilon}}{d{\rho}} - \breve{\lambda}^2 = \breve{\mu}^2 
\quad  ( \breve{\mu} = const  \in \mathbb R )
\end{align*}
 lead to the following system of equations
\begin{gather}
\begin{cases}
\frac{d{^2}{\Xi}}{d{x_0}^2} + \breve{\beta}^2  \Xi = 0, \\[1ex]
x_2^2 \frac{d{^2}{ \Upsilon}}{d{x_2}^2} - \alpha_2 x_2 \frac{d{ \Upsilon}}{d{x_2}}
- (\breve{\lambda}^2  + \breve{\mu}^2)x_2^2 \Upsilon = 0.
\end{cases}
\label{eq-sep-x_2-x_0}
\end{gather}
The first equation of the system~\eqref{eq-sep-x_2-x_0} may be solved using trigonometric functions:
$
\quad \Xi_{\breve{\mu}}(x_0) = b^1_{\breve{\mu}} \cos{\breve{\mu} x_0} +  b^2_{\breve{\mu}} \sin{\breve{\mu} x_0},
$
where $\breve{\mu}\in \mathbb Z$.
 
The second equation of the system~\eqref{eq-sep-x_2-x_0} may be solved using linear independent solutions 
(see, e.g., \cite[Chapter 14, p.535 items 132,126,127]{PolZait:Ordin-2018}):
$$
\Upsilon_{ \breve{\lambda}, \breve{\mu}}(x_2)=
{x_2}^\frac{\alpha_2+1}{2} \left[ a^1_{\breve{\lambda}, \breve{\mu}} I_{\frac{\alpha_2+1}{2}}( \breve{\nu}x_2) +
a^2_{\breve{\lambda}, \breve{\mu}} K_{\frac{\alpha_2+1}{2}}( \breve{\nu}x_2) \right],
$$
 where $I_{\frac{\alpha_2+1}{2}}( \breve{\nu}x_2)$ and  $K_{\frac{\alpha_2+1}{2}}( \breve{\nu}x_2)$
are the modified Bessel functions of the first and second kind of real order ${\frac{\alpha_2 + 1}{2}}$;
$\ \breve{\nu} = \sqrt{ \breve{\lambda}^2 + \breve{\mu}^2}$; 
$\ a^1_{\breve{\lambda}, \breve{\mu}}, a^2_{\breve{\lambda}, \breve{\mu}}= const \in \mathbb R$ (see, e.g., \cite{Watson:1944,Koren:2002}).
 \end{proof}

 \begin{remark}
The Dirichlet problem in a bounded rectangular parallelepiped for eqn~\eqref{alpha_1,2-bihyperbolic-3}
 under the conditions $\alpha_1>0$, $\alpha_2>0$  was studied by Urinov and Karimov in 2023 in a three-dimensional setting \cite{UriKar:2023}. 
Different boundary value problems for elliptic equations with singular coefficients
(see, e.g., \cite{KarimETNieto:2011,NietoKarimET:2012,UrinovKarimovKT:2019,UrinovKarimovKT:2020}) 
may have rich applications in mathematical physics and continuum mechanics.
Two-dimensional analytic models of potential meridional and transverse fields are of particular interest.
\end{remark}

 When $\alpha_1=0$, $\alpha_2 \neq 0$, the equation~\eqref{alpha_1,2-bihyperbolic-3} 
leads to the Weinstein equation in $\mathbb R^3$ (see, e.g., \cite{Leut:CV20,ErOrel:2014}) 
\begin{gather}
 x_2 \Delta{h} - \alpha_2 \frac{\partial{h}}{\partial{x_2}} =0.
\label{alpha-hyperbolic-3}
\end{gather}
 Surprising analytic properties of exact solutions of eqn~\eqref{alpha-hyperbolic-3} have been studied by Leutwiler, Eriksson and Orelma  
in the context of \emph{Hyperbolic function theory in $\mathbb R^3$} (see, e.g., \cite{ErOrel:2014}),
 and later in the context of the theory of \emph{Modified harmonic functions in $\mathbb R^3$} (see, e.g., \cite{Leut:2017-AACA,Leut:2017-CAOT,Leut:2021-MMAS}). 

\begin{definition}
Every exact solution of eqn~\eqref{alpha-hyperbolic-3} under the condition $\alpha_2>0$ 
 in a simply connected open domain $\Lambda \subset \mathbb R^3$ $(x_2 > 0)$ is called $\alpha_2$-hyperbolic harmonic potential in $\Lambda$.
\end{definition}

 When $\alpha_2=1$, fundamentally new analytic properties of exact solutions of eqn~\eqref{alpha-hyperbolic-3} 
have been investigated by Leutwiler and Eriksson-Bique
in the context of \emph{Modified quaternionic analysis in $\mathbb R^3$} (see, e.g., \cite{Leut:CV17,Leut:CV20,Leut:Rud96,ErLe:1998}) 
using the reduced quaternionic power series with complex coefficients. 

 When $\alpha_2 < 0$, exact solutions of~\eqref{alpha-hyperbolic-3} in the context of the theory of \emph{Modified harmonic functions in $\mathbb R^3$}
are referred to as $(-\alpha_2)$-modified harmonic functions (see, e.g., \cite{Leut:2021-MMAS}). 

Let us compare the similarities and differences between eqn~\eqref{eq-axial-hyperbolic-3-alpha} and eqn~\eqref{alpha_1,2-bihyperbolic-3} in Cartesian coordinates. 
This immediately leads to the following formulation.
\begin{proposition} [The first criterion]
  Any  $(\alpha_1, \alpha_2)$-bihyperbolic harmonic potential $h= h(x_0, x_1, x_2)$ in $\Lambda \subset \mathbb R^3$ $(x_1>0, x_2>0)$
 represents an $(\alpha_1+ \alpha_2)$-axial-hyperbolic harmonic potential in $\Lambda$ if and only if 
\begin{gather}
 x_2 \frac{\partial{h}}{\partial{x_1}} = x_1 \frac{\partial{h}}{\partial{x_2}}.
\label{meridional-condition}
\end{gather}
 \end{proposition}

\begin{proof}
 Assume that in~\eqref{eq-axial-hyperbolic-3-alpha} $\alpha = \alpha_1+ \alpha_2$  and $x_1>0$, $x_2>0$.
 Obviously,  $\ x_2 \frac{\partial{h}}{\partial{x_1}} = x_1 \frac{\partial{h}}{\partial{x_2}}$ if and only if
$\ \frac{1}{x_1} \frac{\partial{h}}{\partial{x_1}} = \frac{1}{x_2} \frac{\partial{h}}{\partial{x_2}}$. 

As follows from~\eqref{eq-axial-hyperbolic-3-alpha} and~\eqref{alpha_1,2-bihyperbolic-3},
\begin{gather}
 \Delta{h} = \frac{(\alpha_1+ \alpha_2)x_1}{(x_1^2+x_2^2)} \frac{\partial{h}}{\partial{x_1}} 
 + \frac{(\alpha_1+ \alpha_2) x_2}{(x_1^2+x_2^2)} \frac{\partial{h}}{\partial{x_2}}  =
 \frac{\alpha_1}{x_1} \frac{\partial{h}}{\partial{x_1}} + \frac{\alpha_2}{x_2} \frac{\partial{h}}{\partial{x_2}}.
  \label{Rel-axial-hyperbolic-bihyperbolic-3}
\end{gather}
 Relations~\eqref{Rel-axial-hyperbolic-bihyperbolic-3} imply that
\begin{gather}
  \frac{(\alpha_1+ \alpha_2)x_1^2 - \alpha_1(x_1^2+x_2^2)}{(x_1^2+x_2^2)} 
\frac{1}{x_1} \frac{\partial{h}}{\partial{x_1}}
 = \frac{\alpha_2(x_1^2+x_2^2) - (\alpha_1+ \alpha_2) x_2^2}{(x_1^2+x_2^2)} 
\frac{1}{x_2} \frac{\partial{h}}{\partial{x_2}}.
  \label{alpha-axial-hyperbolic-bihyperbolic-3}
  \end{gather}
Eqn~\eqref{alpha-axial-hyperbolic-bihyperbolic-3} is satisfied if and only if 
the condition~\eqref{meridional-condition} is satisfied.
 \end{proof}

Now let us compare the similarities and differences between
eqns~\eqref{eq-axial-hyperbolic-3-alpha} and~\eqref{alpha_1,2-bihyperbolic-3} in cylindrical coordinates. 
This immediately leads to the following formulation.
 \begin{proposition}
[The second criterion]
Any $(\alpha_1, \alpha_2)$-bihyperbolic harmonic potential $h= h(x_0, x_1, x_2)$ in $\Lambda \subset \mathbb R^3$ $(x_1>0, x_2>0)$
 represents an $(\alpha_1+ \alpha_2)$-axial-hyperbolic harmonic potential in $\Lambda$ if and only if  in cylindrical coordinates 
\begin{gather}
\frac{\partial{h}}{\partial{\theta}} = 0. 
\label{meridional-condition-cyl}
  \end{gather}
\end{proposition}
\begin{proof}
When $\alpha = \alpha_1+ \alpha_2$, eqn~\eqref{eq-axial-hyperbolic-3-alpha} in cylindrical coordinates is written as 
 \begin{gather}
 \rho^2 \left( \frac{\partial{^2}{h}}{\partial{x_0}^2} +  \frac{\partial {^2}{h}}{\partial{\rho}^2} \right) 
- (\alpha_1+ \alpha_2 -1) \rho \frac{\partial{h}}{\partial{\rho}} + \frac{\partial {^2}{h}}{\partial{\theta}^2} = 0.
 \label{eq-axial-hyperbolic-3-alpha-cyl}
\end{gather}
Eqn~\eqref{alpha_1,2-bihyperbolic-3} in cylindrical coordinates is written as
\begin{gather}
\rho^2 \left( \frac{\partial{^2}{h}}{\partial{x_0}^2} +  \frac{\partial {^2}{h}}{\partial{\rho}^2} \right)
  - (\alpha_1 + \alpha_2 -1) \rho \frac{\partial{h}}{\partial{\rho}} + \frac{\partial {^2}{h}}{\partial{\theta}^2}
+ (\alpha_1 \tan{\theta} - \alpha_2 \cot{\theta}) \frac{\partial{h}}{\partial{\theta}} =0.
 \label{alpha_1,2-bihyperbolic-3-cyl}
\end{gather}
This implies that the condition~\eqref{meridional-condition-cyl} is necessary and sufficient.
 \end{proof}

 This approach may be developed in the context of \emph{Hyperbolic function theory  in $\mathbb R^3$},
  the theory of \emph{Modified harmonic functions in $\mathbb R^3$} and \emph{Modified quaternionic analysis in $\mathbb R^3$}
for different values of the parameter $\alpha$.

 \begin{remark} 
 Set of potential functions $h = h(x_0, \rho, \theta)$ in cylindrical coordinates
 within the joint class of $\alpha$-hyperbolic harmonic and $\alpha$-axial-hyperbolic harmonic potentials 
may be replaced by set of scalar functions $g =g(x_0, \rho)$ depending on only two variables. 
\end{remark}

 The general class of  joint exact solutions of eqns~\eqref{eq-axial-hyperbolic-3-alpha-cyl} and~\eqref{alpha_1,2-bihyperbolic-3-cyl}, 
satisfying the condition~\eqref{meridional-condition-cyl}, may be equivalently represented as the general class of exact solutions 
of the elliptic Euler-Poisson-Darboux equation in cylindrical coordinates \cite{Aksenov:2005,Br:Hefei2020}:
\begin{gather}
 \rho \left( \frac{\partial{^2}{g}}{\partial{x_0}^2} +  \frac{\partial {^2}{g}}{\partial{\rho}^2} \right)
  - (\alpha -1) \frac{\partial{g}}{\partial{\rho}}  = 0,
  \label{EPD equation}
  \end{gather}
where, according to \cite{Br:Hefei2020}, $h(x_0, x_1, x_2) := g(x_0, \rho)$, and $\alpha = \alpha_1 + \alpha_2$.

The general class of exact solutions of eqn~\eqref{EPD equation} 
in the context of \emph{GASPT} (see, e.g., \cite{Weinstein:1948-flows,Weinstein:1953,Br:Hefei2020})
is referred to as class of generalized axially symmetric potentials.

The equation
\begin{gather}
g_C(x_0,\rho) = g(x_0,\rho) - C =0
  \label{level-C-axisym}
\end{gather}
allows us to study the specifics of level surfaces of scalar potentials $h$ by means of generalized axially symmetric potentials $g$.

 By the Implicit Function Theorem, equation~\eqref{level-C-axisym} may be solved for $x_0$ in terms of $\rho$ and the constant $C$
in a neighborhood of any point $\vec {\breve{x}} = (\breve{x}_0, \breve{x}_1, \breve{x}_2) \in \Lambda$
at which
 $\frac{\partial{g_C(\breve{x}_0, \breve{\rho})}}{\partial{x_0}} = \frac{\partial{g(\breve{x}_0, \breve{\rho})}}{\partial{x_0}} \neq 0$.
The cylindrical radial variable at $ \vec {\breve{x}} $ is described as $ \breve{\rho} = \sqrt{ {\breve{x}_1}^2+ {\breve{x}_2}^2}$.
Thus, the first derivative of the corresponding function $x_0 = x_0(\rho, C)$ 
may be obtained using implicit differentiation of~\eqref{level-C-axisym}:
\begin{gather*}
\frac{d{x_0}}{d{\rho}} = - \frac{\partial{g}}{\partial{\rho}} \left( \frac{\partial{g}}{\partial{x_0}} \right)^{-1}.
\end{gather*}

Similarly,~\eqref{level-C-axisym} may be solved for $\rho$ in terms of $x_0$ and the constant $C$
in a neighborhood of any point $\vec {\breve{x}} = (\breve{x}_0, \breve{x}_1, \breve{x}_2) \in \Lambda$
at which
 $\frac{\partial{g_C(\breve{x}_0, \breve{\rho})}}{\partial{\rho}} = $
 $\frac{\partial{g(\breve{x}_0, \breve{\rho})}}{\partial{\rho}} \neq 0$.
The first derivative of the corresponding function $\rho = \rho(x_0, C)$  is described as
\begin{gather*}
\frac{d{\rho}}{d{x_0}} = - \frac{\partial{g}}{\partial{x_0}} \left( \frac{\partial{g}}{\partial{\rho}} \right)^{-1}.
\end{gather*}

The first special class of exact solutions of~\eqref{EPD equation}  
under the assumption that  $g(x_0,  \rho) = \Xi(x_0) \Upsilon(\rho)$ may be obtained  
 using hyperbolic functions, power functions and the Bessel functions of the first and second kind \cite{Br:Hefei2020}: 
\begin{align*}
& g(x_0,  \rho) =
\left[ b^1_{\breve{\beta}} \cosh(\breve{\beta} x_0)
+ b^2_{\breve{\beta}} \sinh(\breve{\beta}x_0) \right]
{\rho}^{\frac{\alpha}{2}}
\left[ a^1_{\breve{\beta}} J_{\frac{\alpha}{2}}( \breve{\beta} \rho)
+ a^2_{\breve{\beta}} Y_{\frac{\alpha}{2}}( \breve{\beta} \rho) \right],
\end{align*}
where $\ \breve{\beta}, b^1_{\breve{\beta}}, b^2_{\breve{\beta}},
a^1_{\breve{\beta}}, a^2_{\breve{\beta}} = const \in \mathbb R $.

\begin{proposition}
 The second special class of exact solutions of~\eqref{EPD equation}  
may be obtained using trigonometric functions,  power functions and the modified Bessel functions of the first and second kind:
\begin{align*}
& g(x_0,  \rho) = \sum_{\breve{\gamma}= -\infty}^\infty
\left( b^1_{\breve{\gamma}} \cos(\breve{\gamma} x_0) + b^2_{\breve{\gamma}} \sin(\breve{\gamma} x_0) \right)
{\rho}^{\frac{\alpha}{2}} \left[ a^1_{\breve{\gamma}} I_{\frac{\alpha}{2}}(\breve{\gamma}\rho) +
a^2_{\breve{\gamma}} K_{\frac{\alpha}{2}}(\breve{\gamma}\rho) \right],
\end{align*}
where $\  b^1_{\breve{\gamma}}, b^2_{\breve{\gamma}},  a^1_{\breve{\gamma}}, a^2_{\breve{\gamma}} = const \in \mathbb R $.
\end{proposition}
 \begin{proof}
Consider the second special class of exact solutions of~\eqref{EPD equation}, where $g(x_0,  \rho) = \Xi(x_0) \Upsilon(\rho)$:
$$
\frac{1}{\Xi}  \frac{d{^2}{\Xi}}{d{x_0}^2} +   \frac{1}{ \Upsilon} \frac{d{^2}{ \Upsilon}}{d{\rho}^2} - \frac{(\alpha -1)} { \Upsilon \rho} \frac{d{ \Upsilon}}{d{\rho}} = 0.
$$
Relations
$$
- \frac{1}{\Xi} \frac{d{^2}{\Xi}}{d{x_0}^2} =  \frac{1}{ \Upsilon} \frac{d{^2}{ \Upsilon}}{d{\rho}^2} - \frac{(\alpha -1)} { \Upsilon \rho} \frac{d{ \Upsilon}}{d{\rho}} =  \breve{\gamma}^2 \quad (\breve{\gamma} \in \mathbb R)
$$
 lead to the following system of ordinary differential equations:
\begin{gather}
\begin{cases}
\frac{d{^2}{\Xi}}{d{x_0}^2} + \breve{\gamma}^2  \Xi = 0, \\[1ex]
 \frac{d{^2}{ \Upsilon}}{d{\rho}^2} - \frac{\alpha -1}{\rho} \frac{d{ \Upsilon}}{d{\rho}} - \breve{\gamma}^2 \Upsilon = 0.
\end{cases}
\label{eq-sep-x_2-x_0-hyper-cyl}
\end{gather}
The first equation of~\eqref{eq-sep-x_2-x_0-hyper-cyl} may be solved using  trigonometric  functions:

$ \ \Xi_{\breve{\gamma}}(x_0) = b^1_{\breve{\gamma}} \cos(\breve{\gamma} x_0) + b^2_{\breve{\gamma}} \sin(\breve{\gamma} x_0)$;
 $\ \breve{\gamma}\in \mathbf{Z}$, $\ b^1_{\breve{\gamma}}, b^2_{\breve{\gamma}}= const \in \mathbb R$.

The second equation of~\eqref{eq-sep-x_2-x_0-hyper-cyl} may be solved using linear independent solutions
(see, e.g., \cite[Chapter 14, p.535 items 132,126,127]{PolZait:Ordin-2018}):
$$
\Upsilon_{ \breve{\gamma}}(\rho)=
{\rho}^{\frac{\alpha}{2}} \left[ a^1_{\breve{\gamma}} I_{\frac{\alpha}{2}}(\breve{\gamma}\rho) +
a^2_{\breve{\gamma}} K_{\frac{\alpha}{2}}(\breve{\gamma}\rho) \right],
$$
where  $I_{\frac{\alpha}{2}}(\breve{\gamma}\rho)$ and  $K_{\frac{\alpha}{2}}(\breve{\gamma}\rho)$
are the modified Bessel functions of the first and second kind of real order ${\frac{\alpha}{2}}$;
$\ a^1_{\breve{\gamma}}, a^2_{\breve{\gamma}}= const \in \mathbb R$.
 \end{proof}

Every generalized axially symmetric potential $g = g(x_0, \rho)$ indicates the existence of the Stokes stream function $\hat{g} = \hat{g}(x_0, \rho)$,
which is defined by the generalized Stokes-Beltrami system in the meridian half-plane $(\rho > 0)$
\begin{gather*}
\begin{cases}
  {\rho}^{-(\alpha -1)} \frac{\partial{g}}{\partial{x_0}} = \frac{\partial{\hat{g}}}{\partial{\rho}}, \\[1ex]
  {\rho}^{-(\alpha -1)} \frac{\partial{g}}{\partial{\rho}}=-\frac{\partial{\hat{g}}}{\partial{x_0}}.
  \end{cases}
\end{gather*}
 The Stokes stream function $\hat{g} = \hat{g}(x_0, \rho)$, in contrast to generalized axially symmetric potential, satisfies the following equation:
 \begin{gather}
  \rho \left( \frac{\partial{^2}{\hat{g}}}{\partial{x_0}^2} +  \frac{\partial {^2}{\hat{g}}}{\partial{\rho}^2} \right)
  + (\alpha -1) \frac{\partial{\hat{g}}}{\partial{\rho}} = 0.
\label{Stokes stream}
\end{gather}

\section
{Gradient Systems in $\mathbb R^3$ and $\alpha$-Meridional Mappings of the Second Kind in Continuum Mechanics } 
\label{sec4}

Hereinafter, the vector $\vec V= \mathrm{grad} \ h$ will be identified with a potential velocity field, the scalar potential $h$ with the velocity potential,
the coefficient $\phi$ with the mass density of an inhomogeneous isotropic medium,
and the Jacobian matrix $\mathbf{J}(\vec V)$ with the rate of deformation tensor
 (see, e.g., \cite{LaiRubKr:2010,Reddy:2018,WhiteXue:2021}), respectively.

\begin{definition}
Let $h = h(\vec x)$ be a $C^2$ single-valued velocity potential  
in a simply connected open domain $\Lambda \subset \mathbb R^4=\{(x_0, x_1,x_2)\}$,
 such that $\vec V(\vec x) = \mathrm{grad} \ h(\vec x)$. 
 A smooth gradient system 
\begin{gather}
\frac {d{\vec x}}{dt} = \vec V(\vec x) = \mathrm{grad} \ h(\vec x)
  \label{grad-system-3}
\end{gather} 
  with velocity potential $h$ satisfying the conditions~\eqref{meridional-condition}
is said to be a smooth axisymmetric gradient system in $\Lambda$. 
  \end{definition}

Assume that $ V_l(\vec x) \neq 0$ $(l = 0,1,2)$ in $\Lambda$. 
Then~\eqref{grad-system-3} is written as
  \begin{gather*}
   \frac{dx_0}{V_0} = \frac{dx_1}{V_1} = \frac{dx_2}{V_2} = dt.
  \end{gather*}

Assume that there exists a Jacobi multiplier $M= M(x_0,x_1,x_2)$ for~\eqref{grad-system-3}  in $\Lambda$
when $M$ solves a first-order equation of the form
\begin{gather}
   \frac{\partial{(M V_0)}}{\partial{x_0}} +  \frac{\partial{(M V_1)}}{\partial{x_1}} + \frac{\partial{(M V_2)}}{\partial{x_2}} = 0.
 \label{multiplier-3}
\end{gather}

According to \cite[p.43]{Poincare:1892-1993}, the integral of the order 3
 \begin{gather}
\int_{\Lambda} M(x_0,x_1,x_2)dx_0 dx_1 dx_2 
  \label{integral-3}
 \end{gather}
is an integral invariant of~\eqref{grad-system-3}.
Every $C^1$  function $M= M(x_0,x_1,x_2)$ may be seen as the density of the integral invariant~\eqref{integral-3}
(see, e.g., \cite{Poincare:1892-1993,NemStep:1960,Goriely:2001,BerGiac:2003}). 

According to \cite[p. 63]{Goriely:2001}, the multiplier $M$ may also be seen as the density associated with the invariant measure $\int_{\Lambda} M d \vec x$.

When we deal with gradient systems~\eqref{grad-system-3},
 the first-order equation~\eqref{multiplier-3} is replaced by a second-order equation in divergence form
\begin{gather}
\mathrm{div} \, (M \ \mathrm{grad}{\ h(\vec x)}) =  0.
\label{Liouville-3-Poincare}
\end{gather}

 \begin{remark}
If a coefficient $\phi= \phi(x_0,x_1,x_2)$ within a certain first-order system~\eqref{potential-system-3}
 coincides with a Jacobi multiplier $M= M(x_0,x_1,x_2)$, then eqn~\eqref{Liouville-3} coincides with eqn~\eqref{Liouville-3-Poincare}. 
 \end{remark}

Let us consider basic properties of a smooth gradient system~\eqref{grad-system-3} 
with velocity potential $h$ in the following expanded form:
 \begin{gather}
\begin{cases}
\frac {dx_0}{dt} = V_0(x_0,x_1,x_2) = \frac{\partial{h(x_0,x_1,x_2)}}{\partial{x_0}}, \\[1ex]
\frac {dx_1}{dt} = V_1(x_0,x_1,x_2) = \frac{\partial{h(x_0,x_1,x_2)}}{\partial{x_1}}, \\[1ex]
\frac {dx_2}{dt} = V_2(x_0,x_1,x_2) = \frac{\partial{h(x_0,x_1,x_2)}}{\partial{x_2}}.
\end{cases}
   \label{traject}
 \end{gather}
  These systems in continuum mechanics may be seen as the pathline equations
 (see, e.g., \cite{Sedov:1994,LaiRubKr:2010,Batch:2000,WhiteXue:2021}).
 The system of the streamline equations in continuum mechanics is described as
 \begin{gather}
   \frac{\frac{dx_0}{ds}}{V_0} =  \frac{\frac{dx_1}{ds}}{V_1}  =  \frac{\frac{dx_2}{ds}}{V_2},
   \label{streamline-Acheson}
 \end{gather}
where $s$ is a real-valued parameter.

In general, the systems of equations~\eqref{traject} and~\eqref{streamline-Acheson} are different. 
Nevertheless, the systems~\eqref{traject} and~\eqref{streamline-Acheson} are identical
in the case of the steady state motion of a fluid, where $V_l \neq 0$ $(l = 0,1,2)$ in $\Lambda$.

As noted by White and Xue (\cite{WhiteXue:2021}, p.6),  
``From the point of view of fluid mechanics, all matter consists of only two states, fluid and solid."

 According to \cite[pages 425-426]{NemStep:1960}), every smooth gradient system~\eqref{grad-system-3} may be interpreted 
as a system of equations determining the velocity $\vec V(\vec x)$ of the steady state motion of a fluid filling the space $\mathbb R^3$.
 The integral invariant~\eqref{integral-3} may be interpreted as the mass of the fluid filling a domain $\Lambda$, respectively.

When $V_0 \neq 0$ in $\Lambda$,~\eqref{traject} may be represented as (see, e.g., \cite[pages 43-44]{Sedov:1994})
 \begin{gather}
\begin{cases}
\frac {dx_1}{dx_0} = \frac {V_1(x_0,x_1,x_2)}{V_0(x_0,x_1,x_2)}, \\[1ex]
\frac {dx_2}{dx_0} = \frac {V_2(x_0,x_1,x_2)}{V_0(x_0,x_1,x_2)}.
\end{cases}
 \label{V_0-3}
 \end{gather}

\begin{definition}  
 The set of all points $\vec x = (x_0,x_1,x_2)$, where $V_l(x_0,x_1,x_2; \mu) =0$ $(l = 0,1,2)$  in $\Lambda$,
  is called the $x_l$-nullcline of~\eqref{traject} in $\Lambda$.
  \end{definition}

 According to \cite[p.187]{HirschSmaleDev:2013}, the nullclines may be seen as 
one of the most useful tools for analyzing the behavior of~\eqref{traject} in the context of \emph{Global nonlinear techniques}. 
 In particular,  the intersections of the $x_0$-, $x_1$- and $x_2$-nullclines in $\Lambda$ yield the set of equilibria of~\eqref{traject} in $\Lambda$.

 Fundamentally new properties of potential velocity fields $\vec V$ in inhomogeneous isotropic media 
with the mass density $\phi = \phi(x_0,x_1,x_2)$ may be studied in the context of \emph{Stability theory}
(see, e.g., \cite{Poincare:1885,Chetayev:1990,Gilmore:1993,Kozlov:1993}).

 The sets of equilibria of gradient systems~\eqref{traject} coincide with the critical sets of $h$.
According to \cite[p.9]{Gilmore:1993}, the stability properties of equilibria of~\eqref{traject} 
 are determined from the eigenvalues of the Hessian matrix $\mathbf{H}(h)$.
 When $h$ depends on one or more control parameters, 
the topological and geometric properties of sets of non-Morse critical points of $h$ are of particular interest to 
\emph{Bifurcation theory} and \emph{Catastrophe theory} 
(see, e.g., \cite{Khesin:1990,ArnAfrIlyashShil:1994,Strogatz:2018,Kuznetsov:2023,PosStew,Gilmore:1993}).

Following the concepts expounded by 
Poincar\'{e} \cite{Poincare:1885}, Chetayev \cite{Chetayev:1990}, Gilmore \cite{Gilmore:1993} and Kozlov \cite{Kozlov:1993},
the eigenvalues of the Hessian matrix $\mathbf{H}(h)$ within the critical sets of $h$ may be seen as the 
\emph{Poincar\'{e} stability coefficients} of equilibria of~\eqref{traject}.

\begin{definition} (see, e.g., \cite[p.771]{Kozlov:1993})  Let $\vec x^{**}$ be an equilibrium point of~\eqref{traject}.
The number of positive eigenvalues (counting multiplicities) of the Jacobian matrix $\mathbf{J}(\vec V)$ at $\vec x^{**}$
is called the degree of instability of $\vec x^{**}$.
   \end{definition}

When $\phi = M  = \rho^{-\alpha}$ $(\rho >0)$, eqn~\eqref{EPD equation} leads to
a family of Vekua type systems in the meridian half-plane for different values of $\alpha$  \cite{Br:Hefei2020}:
\begin{gather}
\begin{cases}
 \rho \left( \frac{\partial{u_0}}{\partial{x_0}} - \frac{\partial{u_{\rho}}}{\partial{\rho}} \right) 
+ (\alpha -1) u_{\rho} = 0,
   \\[1ex]
      \frac{\partial{u_0}}{\partial{\rho}}=-\frac{\partial{u_{\rho}}}{\partial{x_0}},
\end{cases}
\label{A_3^alpha system-meridional}
\end{gather} 
where $u_0 = \frac{\partial{g}}{\partial{x_0}}, \quad u_{\rho} = - \frac{\partial{g}}{\partial{\rho}}$.

 The system~\eqref{alpha-axial-hyperbolic-system-3} is reduced to the following two-dimensional system:
\begin{gather}
\begin{cases}
 \rho \left( \frac{\partial{V_0}}{\partial{x_0}} + \frac{\partial{V_{\rho}}}{\partial{\rho}} \right)
 - (\alpha -1) V_{\rho} = 0, 
  \\[1ex]
      \frac{\partial{V_0}}{\partial{\rho}} = \frac{\partial{V_{\rho}}}{\partial{x_0}},
\end{cases}
\label{Bryukhov-vector-meridional}
\end{gather}
where
\begin{gather*}
 V_0= u_0, \quad V_1  = \frac{x_1}{\rho} V_{\rho} = -u_1,  \quad V_2  = \frac{x_2}{\rho} V_{\rho} = -u_2, \quad 
  V_{\rho} = -u_{\rho}.
\end{gather*}

The Jacobian matrix $\mathbf{J}(\vec V)$ of $\vec V = \left(V_0,\frac{x_1}{\rho} V_{\rho},\frac{x_2}{\rho} V_{\rho} \right)$ is expressed as
\begin{gather}
\begin{pmatrix}
 \left[ -\frac{\partial{V_{\rho}}}{\partial{\rho}} +\frac{V_{\rho}}{\rho} (\alpha -1) \right] &  \frac{\partial{V_{\rho}}}{\partial{x_0}} \frac{x_1}{\rho} &
 \frac{\partial{V_{\rho}}}{\partial{x_0}} \frac{x_2}{\rho} \\[1ex]
\frac{\partial{V_{\rho}}}{\partial{x_0}} \frac{x_1}{\rho}  & \left( \frac{\partial{V_{\rho}}}{\partial{\rho}} \frac{x_1^2}{\rho^2}  + \frac{V_{\rho}}{\rho} \frac{x_2^2}{\rho^2}\right)  &
 \left( \frac{\partial{V_{\rho}}}{\partial{\rho}}- \frac{V_{\rho}}{\rho}\right)  \frac{x_1 x_2}{\rho^2}  \\[1ex]
\frac{\partial{V_{\rho}}}{\partial{x_0}} \frac{x_2}{\rho}  & \left( \frac{\partial{V_{\rho}}}{\partial{\rho}}- \frac{V_{\rho}}{\rho}\right)  \frac{x_1 x_2}{\rho^2}  &
\left( \frac{\partial{V_{\rho}}}{\partial{\rho}} \frac{x_2^2}{\rho^2} + \frac{V_{\rho}}{\rho} \frac{x_1^2}{\rho^2}\right)
 \end{pmatrix}
\label{VG tensor-merid}
\end{gather}
The characteristic equation~\eqref{characteristic lambda-3} of~\eqref{VG tensor-merid} is written as   
\begin{gather}
\lambda^3 - \alpha \frac{V_{\rho}}{\rho} \lambda^2 - 
 \left[  \left( \frac{\partial{V_\rho}}{\partial{x_0}}  \right)^2 + \left( \frac{\partial{V_{\rho}}}{\partial{\rho}} \right)^2
- (\alpha -1) \frac{V_{\rho}}{\rho} \left( \frac{\partial{V_{\rho}}}{\partial{\rho}} + \frac{V_{\rho}}{\rho} \right)
 \right]  \lambda \notag \\
 + \frac{V_{\rho}}{\rho} \left[  \left( \frac{\partial{V_\rho}}{\partial{x_0}}  \right)^2 +  \left( \frac{\partial{V_{\rho}}}{\partial{\rho}} \right)^2 
- (\alpha -1) \frac{V_{\rho}}{ \rho} \frac{\partial{V_{\rho}}}{\partial{\rho}} \right] = 0.
 \label{characteristic lambda-alpha}
\end{gather}

\begin{theorem}[see \cite{Br:Hefei2020}]
Roots of~\eqref{characteristic lambda-alpha} are given by the formulas:
\begin{align}
\lambda_{0}
&= \frac{V_{\rho}}{\rho}; \notag\\
\lambda_{1, 2}
&=\frac{(\alpha -1)}{2}  \frac{ V_{\rho}}{ \rho}  \pm \notag\\
&\hspace*{5ex}\sqrt{ \frac{(\alpha -1)^2}{4} \left( \frac{V_{\rho}}{ \rho} \right)^2 - (\alpha -1) \frac{V_{\rho}}{\rho} \frac{\partial{V_{\rho}}}{\partial{\rho}}+
\left( \frac{\partial{V_{\rho}}}{\partial{x_0}}\right)^2 + \left( \frac{\partial{V_{\rho}}}{\partial{\rho}} \right)^2}.
 \label{Roots-alpha}
\end{align}
\end{theorem}
 \begin{remark}
The second formula~\eqref{Roots-alpha} may be simplified: 
\begin{align*}
 \lambda_{1,2}
&= \frac{(\alpha -1)}{2} \frac{V_{\rho}}{\rho}  \pm
 \sqrt{ \left(\frac{\partial{V_{\rho}}}{\partial{x_0}}\right)^2 + 
\left( \frac{\alpha -1}{2} \frac{V_{\rho}}{\rho} - \frac{\partial{V_{\rho}}}{\partial{\rho}} \right)^2}.
\end{align*}
 It implies that the radicand cannot take negative values.
 \end{remark}

As follows from eqn~\eqref{EPD equation} and the first criterion of meridional fields, 
properties of potential meridional fields $\vec V = \mathrm{grad} \ h$
 in cylindrically layered media with the mass density $\phi(\rho) = \rho^{-\alpha}$ 
in $\Lambda$ $(x_1 \neq 0, x_2 \neq 0)$ may be studied in terms of 
gradient systems~\eqref{traject} with $\alpha$-axial-hyperbolic harmonic velocity potential $h$,  
satisfying the condition~\eqref{meridional-condition}. 

The eigenvalues of~\eqref{VG tensor-merid} at any equilibrium point $\vec x^{**}$ of~\eqref{traject}
 with $\alpha$-axial-hyperbolic harmonic velocity potential $h$, satisfying the condition~\eqref{meridional-condition}, 
may be seen as the \emph{Poincar\'{e} stability coefficients} of $\vec x^{**}$.

\begin{theorem}[On the structure of the sets of equilibria of gradient systems]
Assume that the set of equilibria of a gradient system~\eqref{traject} with $\alpha$-axial-hyperbolic harmonic velocity potential $h$,
satisfying the condition~\eqref{meridional-condition}, is not empty in $\Lambda$ $(x_1 \neq 0, x_2 \neq 0)$.
 Then every equilibrium point $\vec x^{**}$ of the system~\eqref{traject} in $\Lambda$ is degenerate.
 The degree of instability of $\vec x^{**}$ is equal to one for any $\alpha$.
\end{theorem}
\begin{proof}
 Set of degenerate points of~\eqref{VG tensor-merid} is provided by two independent equations \cite[p.14]{Br:Hefei2020}:
\begin{align}
{V_{\rho}}=0; \quad
\left(\frac{\partial{V_{\rho}}}{\partial{x_0}}\right)^2 +
\left(\frac{\partial{V_{\rho}}}{\partial{\rho}}\right)^2 -
(\alpha-1)\frac{V_{\rho}}{\rho}\frac{\partial{V_{\rho}}}{\partial{\rho}}=0.
 \label{degenerate-alpha}
\end{align}
Every equilibrium point $\vec x^{**}$ is defined by the condition $\vec V (\vec x^{**}) = 0$. 
As follows from the first equation of~\eqref{degenerate-alpha}, all equilibrium points $\vec x^{**}$
of~\eqref{traject} belong to a set of degenerate points of~\eqref{VG tensor-merid}.
The eigenvalues of~\eqref{VG tensor-merid} at $\vec x^{**}$ are given by the formulas 
\begin{align*}
\lambda_{0}
 &= 0; \notag\\
 \lambda_{1,2} 
 &= \pm \sqrt{ \left(\frac{\partial{V_{\rho}}}{\partial{x_0}}\right)^2 + 
\left( \frac{\partial{V_{\rho}}}{\partial{\rho}} \right)^2}.
\end{align*}
\end{proof}

\begin{definition}
Let $\Lambda \subset \mathbb R^3$  $(x_1 \neq 0, x_2 \neq 0)$ be a simply connected open domain.
 Assume that an exact solution $(u_0, u_1, u_2)$ of the system~\eqref{A_3^alpha-system}, where $\alpha \neq 0$,
satisfies the condition $x_2 u_1 = x_1 u_2$ in $\Lambda$. 
A mapping $u = u_0 + iu_1 + ju_2: \Lambda \rightarrow \mathbb{R}^3$ is called 
$\alpha$-meridional mapping of the first kind, and the corresponding
mapping $ \overline{u} = u_0 - iu_1 - ju_2: \Lambda \rightarrow \mathbb{R}^3$ 
is called $\alpha$-meridional mapping of the second kind.
\end{definition}

\begin{remark}
 Any  $\alpha$-meridional mapping of the second kind may be equivalently represented as a mapping 
$\overline{u} = V_0 + iV_1 + jV_2: \Lambda \rightarrow \mathbb{R}^3$, where $x_2 V_1 = x_1 V_2$.
The Jacobian matrix $\mathbf{J}(\overline{u})$ of $\alpha$-meridional mappings of the second kind
 $\overline{u} = u_0 - iu_1 - ju_2: \Lambda \rightarrow \mathbb{R}^3$ may be identified with 
the Jacobian matrix~\eqref{VG tensor-merid} of potential meridional fields $\vec V$  
in cylindrically layered media with the mass density $\phi( \rho) = \rho^{-\alpha}$.
\end{remark}

 Now let us consider basic properties of harmonic meridional mappings of the second kind
within analytic models of potential meridional fields $\vec V$ in fluid mechanics in the case of a steady flow. 
\begin{definition}
Let $\Lambda \subset \mathbb R^3$  $(x_1 \neq 0, x_2 \neq 0)$ be a simply connected open domain. 
Assume that an exact solution $(u_0, u_1, u_2)$ of the system $(R)$ 
satisfies the condition $x_2 u_1 = x_1 u_2$ in $\Lambda$.
A mapping $u = u_0 + iu_1 + ju_2: \Lambda \rightarrow \mathbb{R}^3$ 
is called harmonic meridional mapping of the first kind, 
and the corresponding mapping $ \overline{u} = u_0 - iu_1 - ju_2: \Lambda \rightarrow \mathbb{R}^3$ 
is called harmonic meridional mapping of the second kind.
\end{definition}

When $\alpha =0$, the system~\eqref{A_3^alpha system-meridional} is expressed as
\begin{gather*}
\begin{cases}
\rho \left( \frac{\partial{u_0}}{\partial{x_0}} - \frac{\partial{u_{\rho}}}{\partial{\rho}} \right) - u_{\rho} = 0,\\[1ex]
\frac{\partial{u_0}}{\partial{\rho}}=-\frac{\partial{u_{\rho}}}{\partial{x_0}}.
\end{cases}
\end{gather*}

The system~\eqref{Bryukhov-vector-meridional} is simplified:
\begin{gather*}
\begin{cases}
\rho \left( \frac{\partial{V_0}}{\partial{x_0}} + \frac{\partial{V_{\rho}}}{\partial{\rho}} \right) + V_{\rho} = 0,  \\[1ex]
\frac{\partial{V_0}}{\partial{\rho}} = \frac{\partial{V_{\rho}}}{\partial{x_0}}.
\end{cases}
\end{gather*}

 Generalized axially symmetric potentials $g = g(x_0, \rho)$ and the Stokes stream functions $\hat{g} = \hat{g}(x_0, \rho)$ satisfy equations
 \begin{gather}
\rho \left( \frac{\partial{^2}{g}}{\partial{x_0}^2} +  \frac{\partial {^2}{g}}{\partial{\rho}^2} \right) + \frac{\partial{g}}{\partial{\rho}} = 0,
\label{EPD equation-0}
\end{gather}
 \begin{gather}
 \rho \left( \frac{\partial{^2}{\hat{g}}}{\partial{x_0}^2} +  \frac{\partial {^2}{\hat{g}}}{\partial{\rho}^2} \right) - \frac{\partial{\hat{g}}}{\partial{\rho}} = 0.
\label{Stokes stream-0}
\end{gather}
 Specific properties of boundary value problems for eqns~\eqref{EPD equation-0} and~\eqref{Stokes stream-0} 
in simply connected domains of the meridian half-plane $(\rho >0)$ and the corresponding applications in fluid mechanics 
have been studied by Plaksa, Shpakivskyi and Gryshchuk in the context of the theory of \emph{Monogenic functions in spaces with commutative multiplication} 
   (see, e.g., \cite{Plaksa:2001,Plaksa:2003,PlakShpak:2023}).

 The geometric properties of the Jacobian matrix~\eqref{VG tensor-merid} of potential meridional fields 
$\vec V$ in homogeneous media
\begin{gather}
\begin{pmatrix}
 \left( -\frac{\partial{V_{\rho}}}{\partial{\rho}} - \frac{V_{\rho}}{\rho} \right) &  \frac{\partial{V_{\rho}}}{\partial{x_0}} \frac{x_1}{\rho} &
 \frac{\partial{V_{\rho}}}{\partial{x_0}} \frac{x_2}{\rho} \\[1ex]
\frac{\partial{V_{\rho}}}{\partial{x_0}} \frac{x_1}{\rho}  &
\left( \frac{\partial{V_{\rho}}}{\partial{\rho}} \frac{x_1^2}{\rho^2} + \frac{E_{\rho}}{\rho} \frac{x_2^2}{\rho^2}\right)  &
 \left( \frac{\partial{V_{\rho}}}{\partial{\rho}}- \frac{V_{\rho}}{\rho}\right)  \frac{x_1 x_2}{\rho^2} \\[1ex]
\frac{\partial{V_{\rho}}}{\partial{x_0}} \frac{x_2}{\rho}  &
\left( \frac{\partial{V_{\rho}}}{\partial{\rho}} - \frac{V_{\rho}}{\rho}\right)  \frac{x_1 x_2}{\rho^2}  &
\left( \frac{\partial{V_{\rho}}}{\partial{\rho}} \frac{x_2^2}{\rho^2} + \frac{V_{\rho}}{\rho} \frac{x_1^2}{\rho^2}\right)
 \end{pmatrix}
\label{VG tensor-merid-0}
\end{gather}
were presented in 2021 \cite{Br:Hefei2020}. According to \cite[pages 210-211]{LavSh-Hydro}, 
analytic models of potential fields with axial symmetry are of particular interest to hydrodynamics problems.

The principal invariants of~\eqref{VG tensor-merid-0} are described as
\begin{align*}
I_{\mathbf{J}(\vec V)} 
&\equiv \mathrm{div} \, \vec V  = 0, \\[1ex]
II_{\mathbf{J}(\vec V)}
&= - \left[ \left( \frac{\partial{V_\rho}}{\partial{x_0}}  \right)^2 +  \left( \frac{\partial{V_{\rho}}}{\partial{\rho}} \right)^2 \right]
 - \frac{V_{\rho}}{\rho} \frac{\partial{V_{\rho}}}{\partial{\rho}} - \left( \frac{V_{\rho}}{ \rho} \right)^2,\\[1ex]
III_{\mathbf{J}(\vec V)} \equiv \det\mathbf{J}(\vec V)
&= - \frac{V_{\rho}}{\rho} \left[  \left( \frac{\partial{V_\rho}}{\partial{x_0}}  \right)^2 +  \left( \frac{\partial{V_{\rho}}}{\partial{\rho}} \right)^2 
 + \frac{V_{\rho}}{ \rho} \frac{\partial{V_{\rho}}}{\partial{\rho}} \right].
\end{align*}

 The third principal invariant of~\eqref{VG tensor-merid-0} satisfies the inequality $III_{\mathbf{J}(\vec V)} < 0$ under the conditions 
\begin{gather*}
 V_{\rho} > 0; \quad
\left( \frac{\partial{V_\rho}}{\partial{x_0}} \right)^2 + \left( \frac{\partial{V_{\rho}}}{\partial{\rho}} \right)^2
 + \frac{V_{\rho}}{ \rho} \frac{\partial{V_{\rho}}}{\partial{\rho}} > 0,
\end{gather*}
as well as  under the conditions 
\begin{gather*}
 V_{\rho} < 0; \quad
\left( \frac{\partial{V_\rho}}{\partial{x_0}} \right)^2 + \left( \frac{\partial{V_{\rho}}}{\partial{\rho}} \right)^2
 + \frac{V_{\rho}}{ \rho} \frac{\partial{V_{\rho}}}{\partial{\rho}} < 0.
\end{gather*}

The characteristic equation~\eqref{characteristic lambda-alpha} is written as incomplete cubic equation
\begin{gather}
\lambda^3 - 
 \left[ \left( \frac{\partial{V_\rho}}{\partial{x_0}}  \right)^2 + \left( \frac{\partial{V_{\rho}}}{\partial{\rho}} \right)^2
 + \frac{V_{\rho}}{\rho} \left( \frac{\partial{V_{\rho}}}{\partial{\rho}} + \frac{V_{\rho}}{\rho} \right)
 \right]  \lambda \notag \\
 + \frac{V_{\rho}}{\rho} \left[  \left( \frac{\partial{V_\rho}}{\partial{x_0}}  \right)^2 +  \left( \frac{\partial{V_{\rho}}}{\partial{\rho}} \right)^2 + \frac{V_{\rho}}{ \rho} \frac{\partial{V_{\rho}}}{\partial{\rho}} \right] = 0.
\label{characteristic-0}
\end{gather}

  Roots of~\eqref{characteristic-0} are  given by the formulas
\begin{gather*}
\lambda_{0} = \frac{V_{\rho}}{\rho}; \quad
\lambda_{1, 2} = - \frac{V_{\rho}}{2 \rho}  \pm  \sqrt{\left( \frac{V_{\rho}}{2 \rho} \right)^2 +
\frac{V_{\rho}}{\rho} \frac{\partial{V_{\rho}}}{\partial{\rho}} + \left( \frac{\partial{V_{\rho}}}{\partial{x_0}}\right)^2 + \left( \frac{\partial{V_{\rho}}}{\partial{\rho}} \right)^2}.
\end{gather*}
Set of degenerate points of the Jacobian matrix~\eqref{VG tensor-merid-0} is provided by two independent equations:
$$
{V_{\rho}}=0; \quad 
\left(\frac{\partial{V_{\rho}}}{\partial{x_0}}\right)^2 +\left(\frac{\partial{V_{\rho}}}{\partial{\rho}}\right)^2 + \frac{V_{\rho}}{\rho}\frac{\partial{V_{\rho}}}{\partial{\rho}}=0.
$$

Potential meridional velocity fields $\vec V$ in fluid mechanics in the case of a steady flow may be  
 seen as irrotational axisymmetric flow fields (see, e.g., \cite{Batch:2000,WhiteXue:2021}).
According to (\cite{WhiteXue:2021}, p.17), ``Foremost among the properties of a flow is the velocity field $\vec V$. 
In fact, determining the velocity is often tantamount to solving a flow problem, since other properties follow directly from the velocity field." 

 Any equilibrium point $\vec x^{**}$ of gradient systems~\eqref{traject} with harmonic velocity potential 
may be seen as stagnation point (see, e.g., \cite {Acheson,WhiteXue:2021}).
 By Theorem 6, the degree of instability of $\vec x^{**}$ is equal to one.

When $\alpha =0$, the first special class of exact solutions of~\eqref{EPD equation}  
under the assumption that  $g(x_0,  \rho) = \Xi(x_0) \Upsilon(\rho)$ may be obtained  
 using hyperbolic functions and the Bessel functions of the first and second kind \cite{Br:Hefei2020}: 
\begin{align*}
& g(x_0,  \rho) =
\left[ b^1_{\breve{\beta}} \cosh(\breve{\beta} x_0)
+ b^2_{\breve{\beta}} \sinh(\breve{\beta}x_0) \right]
\left[ a^1_{\breve{\beta}} J_{0}( \breve{\beta} \rho)
+ a^2_{\breve{\beta}} Y_{0}( \breve{\beta} \rho) \right],
\end{align*}
where $\ \breve{\beta}, b^1_{\breve{\beta}}, b^2_{\breve{\beta}},
a^1_{\breve{\beta}}, a^2_{\breve{\beta}} = const \in \mathbb R $.

The geometric properties of level hypersurfaces of velocity potentials $h$ may be characterized by means of generalized axially symmetric potentials $g$:
\begin{gather*}
b^1_{\breve{\beta}} \cosh(\breve{\beta} x_0) + b^2_{\breve{\beta}} \sinh(\breve{\beta}x_0) = 
\frac{C}{a^1_{\breve{\beta}} J_{0}( \breve{\beta} \rho) + a^2_{\breve{\beta}} Y_{0}( \breve{\beta} \rho)}.
\end{gather*}

 \begin{example}
Asume that $ b^1_{\breve{\beta}} = b^2_{\breve{\beta}} =1, a^1_{\breve{\beta}} =1, a^2_{\breve{\beta}} = 0$.
The essential properties of the corresponding axisymmetric gradient system with harmonic velocity potential in $\Lambda$ $(\rho >0)$ 
may be demonstrated explicitly using the Bessel function of the first kind $J_{0}(\breve{\beta} \rho)$:
\begin{gather*}
g(x_0, \rho) = e^{\breve{\beta} x_0} \frac{\cos ( \breve{\beta} \rho)}{\rho}.
\end{gather*}

\begin{gather*}
V_0 = \frac{\partial{g}}{\partial{x_0}} = \breve{\beta}e^{\breve{\beta} x_0} J_{0}( \breve{\beta} \rho),
\quad V_{\rho} = \frac{\partial{g}}{\partial{\rho}} =
- \breve{\beta} e^{\breve{\beta} x_0} J_{1}(\breve{\beta} \rho).
\end{gather*}

Potential meridional fields $\vec V$ are represented as
\begin{gather*}
\vec V =(V_0,\frac{x_1}{\rho}V_{\rho},\frac{x_2}{\rho} V_{\rho})=
 \breve{\beta} e^{\breve{\beta} x_0}\left( J_{0}( \breve{\beta} \rho), -\frac{x_1}{\rho} J_{1}( \breve{\beta} \rho), -\frac{x_2}{\rho} J_{1}(\breve{\beta} \rho) \right).
\end{gather*}

When $J_{0}( \breve{\beta} \rho) \neq 0$ in $\Lambda$,~\eqref{V_0-3} is represented as   
\begin{gather*}
\begin{cases}
\frac {dx_1}{dx_0} = - \frac { x_1 J_{1}( \breve{\beta} \rho)}{ \rho J_{0}( \breve{\beta} \rho)}, \\[1ex]
\frac {dx_2}{dx_0} = - \frac { x_2 J_{1}( \breve{\beta} \rho)}{ \rho J_{0}( \breve{\beta} \rho)}.
\end{cases}
 \end{gather*}

Relations  
\begin{gather*}
\begin{cases}
   J_{0}( \breve{\beta} \rho) = 0,  \\[1ex]
  J_{1}( \breve{\beta} \rho) =0
\end{cases}
\end{gather*}
imply that zero set of $\vec V$ is empty (see, e.g., \cite{Watson:1944}).

\begin{gather*}
\frac{\partial{V_{\rho}}}{\partial{x_0}}=
 - \frac{\breve{\beta}^2}{2} e^{\breve{\beta}x_0}[J_{0}( \breve{\beta} \rho) - J_{2}( \breve{\beta} \rho)], \\
 \frac{\partial{V_{\rho}}}{\partial{\rho}} = -\frac{\partial{V_0}}{\partial{x_0}} - \frac{V_{\rho}}{\rho}
 = -\breve{\beta}^2 e^{\breve{\beta} x_0} J_{0}( \breve{\beta} \rho)
 + \breve{\beta} e^{\breve{\beta}x_0} \frac{ J_{1}(\breve{\beta} \rho) }{\rho}.
\end{gather*}

The eigenvalues of~\eqref{VG tensor-merid-0} are  given by the formulas
\begin{align*}
\lambda_{0}
&= -\breve{\beta} e^{\breve{\beta} x_0} \frac{J_{1}(\breve{\beta} \rho) }{\rho}; \notag\\
\lambda_{1, 2}
&= -\breve{\beta} e^{\breve{\beta} x_0} \frac{J_{1}(\breve{\beta} \rho) }{2\rho} \pm \notag\\
&\hspace*{5ex}  \breve{\beta} e^{\breve{\beta} x_0} \sqrt{
\left( \breve{\beta}^2 + \frac{1}{4\rho^2} \right) J_{1}^2(\breve{\beta} \rho)
 + \frac{\breve{\beta} J_{1}(\breve{\beta} \rho) }{2\rho}
 \left[J_{0}(\breve{\beta} \rho) - J_{2}(\breve{\beta} \rho)\right]
 +  \frac{\breve{\beta}^2}{4}\left[J_{0}(\breve{\beta}\rho) - J_{2}(\breve{\beta}\rho) \right]^2}.
\end{align*}

 Set of degenerate points of~\eqref{VG tensor-merid-0} is provided by two independent equations:
\begin{gather*}
  J_{1}(\breve{\beta} \rho)  =0; \\[1ex]
 \breve{\beta} J_{1}^2(\breve{\beta} \rho)
 + \frac{ J_{1}(\breve{\beta} \rho) }{2\rho}
 \left[J_{0}(\breve{\beta} \rho) - J_{2}(\breve{\beta} \rho)\right]
 +  \frac{\breve{\beta}}{4}\left[J_{0}(\breve{\beta}\rho) - J_{2}(\breve{\beta}\rho) \right]^2 =0.
\end{gather*}

 \end{example}

\section
 {$1$-Meridional Mappings of the Second Kind and the Radially Holomorphic Potential in $\mathbb R^3$}
\label{sec5}
Consider basic properties of $1$-meridional mappings of the second kind 
within analytic models of potential meridional fields $\vec V$ in cylindrically layered media with the mass density $\phi( \rho) = \rho^{-1}$ .

 When $\alpha =1$, the system~\eqref{A_3^alpha system-meridional} 
 leads the Cauchy-Riemann type system in the meridian half-plane (see, e.g., \cite{Br:Hefei2020})
\begin{gather}
\begin{cases}
\frac{\partial{u_0}}{\partial{x_0}} - \frac{\partial{u_{\rho}}}{\partial{\rho}}= 0,\\[1ex]
\frac{\partial{u_0}}{\partial{\rho}}=-\frac{\partial{u_{\rho}}}{\partial{x_0}}.
\end{cases}
\label{CR-merid}
\end{gather}

 The system~\eqref{Bryukhov-vector-meridional} is simplified:
\begin{gather}
\begin{cases}
\frac{\partial{V_0}}{\partial{x_0}} + \frac{\partial{V_{\rho}}}{\partial{\rho}}  = 0,\\[1ex]
\frac{\partial{V_0}}{\partial{\rho}} = \frac{\partial{V_{\rho}}}{\partial{x_0}}.
\end{cases}
\label{Bryukhov-vector-1}
\end{gather}

The Jacobian matrix~\eqref{VG tensor-merid} is expressed as
\begin{gather}
\begin{pmatrix}
  -\frac{\partial{V_{\rho}}}{\partial{\rho}} &  \frac{\partial{V_{\rho}}}{\partial{x_0}} \frac{x_1}{\rho} &
 \frac{\partial{V_{\rho}}}{\partial{x_0}} \frac{x_2}{\rho} \\[1ex]
\frac{\partial{V_{\rho}}}{\partial{x_0}} \frac{x_1}{\rho}  & \left( \frac{\partial{V_{\rho}}}{\partial{\rho}} \frac{x_1^2}{\rho^2}  + \frac{V_{\rho}}{\rho} \frac{x_2^2}{\rho^2}\right)  &
 \left( \frac{\partial{V_{\rho}}}{\partial{\rho}}- \frac{V_{\rho}}{\rho}\right)  \frac{x_1 x_2}{\rho^2}  \\[1ex]
\frac{\partial{V_{\rho}}}{\partial{x_0}} \frac{x_2}{\rho}  & \left( \frac{\partial{V_{\rho}}}{\partial{\rho}}- \frac{V_{\rho}}{\rho}\right)  \frac{x_1 x_2}{\rho^2}  &
\left( \frac{\partial{V_{\rho}}}{\partial{\rho}} \frac{x_2^2}{\rho^2} + \frac{V_{\rho}}{\rho} \frac{x_1^2}{\rho^2}\right)
 \end{pmatrix}
\label{VG tensor-merid-1}
\end{gather}

The characteristic equation~\eqref{characteristic lambda-alpha} of~\eqref{VG tensor-merid-1} is written as  
\begin{gather}
\lambda^3 - \frac{V_{\rho}}{\rho} \lambda^2 - \left[ \left( \frac{\partial{V_\rho}}{\partial{x_0}}  \right)^2 
+ \left( \frac{\partial{V_{\rho}}}{\partial{\rho}} \right)^2 \right]  \lambda
 + \frac{V_{\rho}}{\rho} \left[ \left( \frac{\partial{V_\rho}}{\partial{x_0}} \right)^2 + \left( \frac{\partial{V_{\rho}}}{\partial{\rho}} \right)^2 \right] = 0.
\label{characteristic-1}
\end{gather}

 Roots of~\eqref{characteristic-1} are  given by the formulas
\begin{gather*}
\lambda_{0} = \frac{V_{\rho}}{\rho}; \quad
\lambda_{1, 2} = \pm  \sqrt{ \left(\frac{\partial{V_{\rho}}}{\partial{x_0}}\right)^2 +
\left(\frac{\partial{V_{\rho}}}{\partial{\rho}}\right)^2}.
\end{gather*}
Set of degenerate points of~\eqref{VG tensor-merid-1} is provided by two independent equations \cite{Br:Hefei2020}:
\begin{gather}
{V_{\rho}}=0; \quad
\left(\frac{\partial{V_{\rho}}}{\partial{x_0}}\right)^2 +
\left(\frac{\partial{V_{\rho}}}{\partial{\rho}}\right)^2 =0.
\label{degenerate-1}
\end{gather}

Let us turn our attention to the concept of \emph{Radially holomorphic functions} 
introduced by G\"{u}rlebeck, Habetha and Spr\"{o}ssig in 2008 \cite{GuHaSp:2008}.
\begin{definition}
The radial differential operator in $\mathbb R^3$ is defined by formula
\begin{gather*}
\partial_{rad}{G}  := \frac{1}{2} \left( \frac{\partial{}}{\partial{x_0}}- I \frac{\partial{}}{\partial{\rho}} \right) G := G'  \quad (G = g + I \hat{g}).
\end{gather*}
Every reduced quaternion-valued function $G = g + I \hat{g}$ satisfying 
a Cauchy-Riemann type differential equation in simply connected open domains 
$\Lambda \subset \mathbb R^3$  $(x_1 \neq 0, x_2 \neq 0)$
\begin{gather}
\overline{\partial}_{rad}{G} :=  \frac{1}{2} \left( \frac{\partial{}}{\partial{x_0}} 
+ I \frac{\partial{}}{\partial{\rho}} \right) G = 0
\label{CRDO-con}
\end{gather}
is called radially holomorphic in $\Lambda$. The reduced quaternion-conjugate function $\overline{G} = g - I \hat{g}$
 is called radially anti-holomorphic in $\Lambda$.
\end{definition}

As noted in \cite{GuHaSp:2008,Br:Hefei2020}, the basic elementary radially holomorphic functions in $\mathbb R^3$ 
may be viewed as elementary functions of the reduced quaternionic variable
 $ G = G(x) = g(x_0, \rho) + I \hat{g}_{\rho}(x_0, \rho)$ satisfying relations
\begin{gather*}
 [G(x) = x^{n} := r^{n}(\cos{n\varphi} + I \sin{n\varphi})]' = F(x) = n x^{n-1};\\[1ex]
 [G(x) = e^{x} := e^{x_0}(\cos{\rho}+I \sin{\rho})]'= F(x) = e^{x}; \\[1ex]
 [G(x) = \cos{x} := \frac{1}{2}( e^{Ix} + e^{-Ix})]'= F(x) = - \sin{x};\\[1ex]
 [G(x) = \sin{x} := - \frac{I}{2}( e^{Ix} - e^{-Ix})]'= F(x) = \cos{x}; \\[1ex]
 [G(x) = \ln{x} := \ln r +I \varphi]' = F(x) = x^{-1}.
\end{gather*}

\begin{definition}
Let $G = g + I \hat{g}$ be a radially holomorphic function in $\Lambda$ satisfying the following differential equation:
\begin{gather}
 G' = F,
\label{primitive}
\end{gather}
where $F = u_0 + I u_{\rho}$ is a radially holomorphic function
with components $u_0$, $u_{\rho}$ satisfying the system~\eqref{CR-merid} in $\Lambda$. 
The function $G$ is called radially holomorphic primitive of $F$ in $\Lambda$.
\end{definition}

The radially anti-holomorphic function $ \overline{F(x)} = u_0(x_0, \rho) - I u_{\rho}(x_0, \rho)$  
 may be equivalently represented as $\overline{F(x)} = V_0(x_0, \rho) + I V_{\rho}(x_0, \rho)$,
where the components $V_0$ and $V_{\rho}$ satisfy the system~\eqref{Bryukhov-vector-1}.

\begin{definition}
Radially holomorphic primitive $G = g + I \hat{g}$ in $\Lambda$ into the framework of the system~\eqref{Bryukhov-vector-1} 
is said to be the radially holomorphic potential in $\Lambda$.
\end{definition}

New tools of the radially holomorphic potential in $\mathbb R^3$ 
were developed by the author in 2021 \cite{Br:Hefei2020} as the radially holomorphic extension of 
 tools of the complex potential in the context of \emph{GASPT}.
According to \cite{Br:Hefei2020}, general class of exact solutions of eqn~\eqref{CRDO-con} 
is equivalently represented as general class of exact solutions 
of the Cauchy-Riemann type system in the meridian half-plane~\eqref{CR-merid} concerning functions
$g = g(x_0, \rho)$, $\hat{g} = \hat{g}(x_0, \rho)$:
\begin{gather}
\begin{cases}
\frac{\partial{g}}{\partial{x_0}}-\frac{\partial{\hat{g}}}{\partial{\rho}}=0,\\[1ex]
\frac{\partial{g}}{\partial{\rho}}=-\frac{\partial{\hat{g}}}{\partial{x_0}}.
\end{cases}
\label{invar}
\end{gather}

\begin{proposition}
Let  $G = g + I \hat{g}$ be a radially holomorphic function in $\Lambda$.
Then  any reduced quaternion-valued expression of the form 
$\Theta = \vartheta + I \hat{\vartheta} = I G = -\hat{g} + I g $ 
  yields also a radially holomorphic function in $\Lambda$.
\end{proposition}
 \begin{proof} 
The system~\eqref{invar} may be equivalently represented as 
\begin{gather*}
\begin{cases}
 \frac{\partial{(-\hat{g})}}{\partial{\rho}} = - \frac{\partial{g}}{\partial{x_0}},\\[1ex]
 \frac{\partial{(-\hat{g})}}{\partial{x_0}} -  \frac{\partial{g}}{\partial{\rho}} =0.
\end{cases}
\end{gather*}
Mutual permutation of equations leads to the Cauchy-Riemann type system
\begin{gather*}
\begin{cases}
\frac{\partial{(-\hat{g})}}{\partial{x_0}} -  \frac{\partial{g}}{\partial{\rho}} =0, \\[1ex] 
 \frac{\partial{(-\hat{g})}}{\partial{\rho}} = - \frac{\partial{g}}{\partial{x_0}},
\end{cases}
\end{gather*}
where $\vartheta = -\hat{g}$, $\hat{\vartheta} = g$. 
  \end{proof}

 As may be seen, generalized axially symmetric potential $g = g(x_0, \rho)$ 
and the Stokes stream function $\hat{g} = \hat{g}(x_0, \rho)$ within analytic models of potential meridional fields 
$\vec V$ in $\mathbb R^3=\{(x_0, x_1,x_2)\}$ may reverse roles in the context of \emph{GASPT},
similar to the velocity potential and the stream function in the context of \emph{Complex analysis} 
(see, e.g., \cite{KochinKibelRoze:1964,LavSh:1987,Batch:2000,Acheson,WhiteXue:2021}).
 If $g = g(x_0, \rho)$  and $\hat{g} = \hat{g}(x_0, \rho)$ both exist in $\Lambda$,
then the streamlines and equipotential lines are everywhere mutually perpendicular except at zeros of $\vec V$ in $\Lambda$. 

\begin{corollary}  
A linear superposition of a radially holomorphic function $G$ and the radially holomorphic function $G$ with coefficient $I$ 
 yields a radially holomorphic potential in the context of \emph{GASPT}.
\end{corollary}
\begin{proof} 
The system of partial differential equations~\eqref{invar} is linear.
Therefore, any sum of several solutions of~\eqref{invar} is also a solution. 
  \end{proof}

 The main properties of axisymmetric gradient systems~\eqref{traject} with $1$-axial-hyperbolic harmonic velocity potential  
 may be described in terms of a special class of reduced-quaternion-valued ordinary differential equations 
\begin{gather}
  \frac {dx}{dt}  = \overline{F(x)},
 \label{Petr-4}
  \end{gather}
where  $ \overline{F(x)} = u_0 - i\frac{x_1}{\rho}u_{\rho} - j\frac{x_2}{\rho}u_{\rho}  
= u_0 - I u_{\rho}  = V_0 + I V_{\rho}$.

Zero sets of radially anti-holomorphic functions $\overline{F(x)} = V_0 + I V_{\rho}$
coincide with the sets of equilibria of~\eqref{traject} with $1$-axial-hyperbolic harmonic velocity potentials.

 Set of velocity potentials $h = h(x_0, \rho, \theta)$ in cylindrical coordinates under the conditions of~\eqref{meridional-condition-cyl}
may be replaced by a set of generalized axially symmetric potentials $g =g(x_0, \rho)$, where 
\begin{gather*}
\frac{\partial{^2}{g}}{\partial{x_0}^2} +  \frac{\partial{^2}{g}}{\partial{\rho}^2} = 0. 
\end{gather*}

The first special class of generalized axially symmetric potentials $g =g(x_0, \rho)$ 
  may be obtained using trigonometric functions,  power functions and the Bessel functions of the first and second kind of order ${\frac{1}{2}}$: 
\begin{gather*}
 g(x_0,  \rho) = \Xi(x_0)  \Upsilon(\rho) = \notag \\
\left[ b^1_{\breve{\beta}} \cosh(\breve{\beta} x_0)
+ b^2_{\breve{\beta}} \sinh(\breve{\beta}x_0) \right]
{\rho}^{\frac{1}{2}}
\left[ a^1_{\breve{\beta}} J_{\frac{1}{2}}( \breve{\beta} \rho)
+ a^2_{\breve{\beta}} Y_{\frac{1}{2}}( \breve{\beta} \rho) \right], 
\end{gather*}
where $\ \breve{\beta}, b^1_{\breve{\beta}}, b^2_{\breve{\beta}},
a^1_{\breve{\beta}}, a^2_{\breve{\beta}} = const \in \mathbb R $.

 Assume that 
\begin{gather}
\begin{cases}
  \Xi_{\breve{\beta}}(x_0) = b^1_{\breve{\beta}} \cosh(\breve{\beta} x_0)
+ b^2_{\breve{\beta}} \sinh(\breve{\beta}x_0);
\quad \breve{\beta}, b^1_{\breve{\beta}}, b^2_{\breve{\beta}}= const \in \mathbb R, \\[1ex]
  \Upsilon_{\breve{\beta}}(\rho) = {\rho}^\frac{1}{2}
\left[ a^1_{\breve{\beta}} J_{\frac{1}{2}}( \breve{\beta} \rho)
+ a^2_{\breve{\beta}} Y_{\frac{1}{2}}( \breve{\beta} \rho) \right];
  \quad a^1_{\breve{\beta}}$, $a^2_{\breve{\beta}}= const \in \mathbb R,
 \end{cases}
 \label{EPD 1}
\end{gather}
where $\breve{\beta}> 0$, $b^1_{\breve{\beta}} = b^2_{\breve{\beta}} =1$,
 $\ a^1_{\breve{\beta}} = -A^1_{\breve{\beta}} \sqrt{\frac{\pi \breve{\beta}}{2}}$,
 $\ a^2_{\breve{\beta}} = -A^2_{\breve{\beta}} \sqrt{\frac{\pi \breve{\beta}}{2}}$;
$\quad A^1_{\breve{\beta}}, A^2_{\breve{\beta}} = const \in \mathbb R.$

The Bessel functions of the first and second kind of half-integer order may be expressed
in finite terms using algebraic and trigonometric functions (see, e.g., \cite[pages 52-55]{Watson:1944} and \cite[pages 15-17]{Koren:2002}).

$ \ J_{\frac{1}{2}}(\breve{\beta} \rho)
= \sqrt{\frac{2}{\pi \breve{\beta} }} \rho^{-\frac{1}{2}} \sin ( \breve{\beta} \rho)$ and
$ \ Y_{\frac{1}{2}}(\breve{\beta} \rho)
= -\sqrt{\frac{2}{\pi \breve{\beta} }} \rho^{-\frac{1}{2}} \cos ( \breve{\beta} \rho)$.

 Formulas~\eqref{EPD 1} imply that
\begin{gather*}
\begin{cases}
  \Xi_{\breve{\beta}}(x_0) = e^{\breve{\beta}x_0}, \\[1ex]
  \Upsilon_{\breve{\beta}}(\rho) = - A^1_{\breve{\beta}} \sin ( \breve{\beta} \rho) + A^2_{\breve{\beta}} \cos ( \breve{\beta} \rho),
 \end{cases}
\end{gather*}
 and set of generalized axially symmetric potentials is described as
\begin{gather}
 g(x_0,  \rho) =  e^{\breve{\beta} x_0} \left[ -A^1_{\breve{\beta}} \sin ( \breve{\beta} \rho)
+ A^2_{\breve{\beta}} \cos ( \breve{\beta} \rho) \right].
  \label{GASP-1}
\end{gather}

The geometric properties of level surfaces of $h$  may be characterized by means of generalized axially symmetric potentials $g$:
\begin{gather*}
e^{-\breve{\beta} x_0} = \frac{- A^1_{\breve{\beta}} \sin ( \breve{\beta} \rho) + A^2_{\breve{\beta}} \cos ( \breve{\beta} \rho)}{C}.
\end{gather*}

\begin{example}

Let us consider the radially holomorphic potentials $G = G(x)$  represented by a superposition of
the radially holomorphic exponential function $e^{\breve{\beta} x}$ and function $e^{\breve{\beta} x}$ with coefficient $I$:
\begin{align*}
  G(x) = A^2_{\breve{\beta}} e^{\breve{\beta}x} + A^1_{\breve{\beta}} I e^{\breve{\beta}x}
= A^2_{\breve{\beta}} e^{\breve{\beta}x_0}[\cos(\breve{\beta}\rho)+I \sin(\breve{\beta}\rho)]
 + A^1_{\breve{\beta}} e^{\breve{\beta}x_0}[-\sin(\breve{\beta}\rho)+I \cos(\breve{\beta}\rho)].
\end{align*}

 It is easy to check that 
 \begin{gather}
G(x) = g + I \hat{g} =  e^{\breve{\beta} x_0} \left[- A^1_{\breve{\beta}} \sin ( \breve{\beta} \rho)
+ A^2_{\breve{\beta}} \cos ( \breve{\beta} \rho) \right] +
 I e^{\breve{\beta} x_0}  \left[ A^1_{\breve{\beta}} \cos ( \breve{\beta} \rho)
+ A^2_{\breve{\beta}} \sin ( \breve{\beta} \rho) \right].
 \label{exponential}
\end{gather}

In this case, the scalar parts of the radially holomorphic potentials~\eqref{exponential}
coincide with generalized axially symmetric potentials of the form~\eqref{GASP-1}.

The radially holomorphic functions $F = F(x)$ are described as
\begin{gather*}
 F(x) = \breve{\beta} \left(A^2_{\breve{\beta}} e^{\breve{\beta}x} + A^1_{\breve{\beta}} I e^{\breve{\beta}x} \right).
\end{gather*}
\begin{gather*} 
\overline{F(x)} = \breve{\beta} e^{\breve{\beta} x_0}  \left[- A^1_{\breve{\beta}} \sin ( \breve{\beta} \rho)
+ A^2_{\breve{\beta}} \cos ( \breve{\beta} \rho) \right]
- I \breve{\beta} e^{\breve{\beta} x_0}  \left[ A^1_{\breve{\beta}} \cos ( \breve{\beta} \rho)
+ A^2_{\breve{\beta}} \sin ( \breve{\beta} \rho) \right] .
\end{gather*}

The essential properties of the corresponding axisymmetric gradient systems with $1$-axial-hyperbolic harmonic velocity potential  
in $\Lambda$ $(\rho >0)$ may be demonstrated explicitly using set of generalized axially symmetric potentials of the form~\eqref{GASP-1}.

\begin{gather*}
V_0 = \breve{\beta} e^{\breve{\beta} x_0} \left[- A^1_{\breve{\beta}} \sin ( \breve{\beta} \rho) + A^2_{\breve{\beta}} \cos ( \breve{\beta} \rho) \right], \qquad
 V_{\rho} =  - \breve{\beta} e^{\breve{\beta} x_0} \left[ A^1_{\breve{\beta}} \cos ( \breve{\beta} \rho)
+ A^2_{\breve{\beta}} \sin ( \breve{\beta} \rho) \right]. 
\end{gather*}

Equation~\eqref{Petr-4} is written as 
\begin{gather*}
  \frac {dx}{dt}  
= \breve{\beta}  e^{\breve{\beta} x_0} 
\left\{ \left[- A^1_{\breve{\beta}} \sin ( \breve{\beta} \rho) + A^2_{\breve{\beta}} \cos ( \breve{\beta} \rho) \right] 
-  I \left[ A^1_{\breve{\beta}} \cos ( \breve{\beta} \rho) + A^2_{\breve{\beta}} \sin ( \breve{\beta} \rho) \right] \right\}.
  \end{gather*}

When $ \sin ( \breve{\beta} \rho) \neq 0$ in $\Lambda$,~\eqref{V_0-3} is represented as  
\begin{gather*}
\begin{cases}
\frac {dx_1}{dx_0} = \frac { x_1 [ A^1_{\breve{\beta}} \cot ( \breve{\beta} \rho) + A^2_{\breve{\beta}} ]}
{ \rho \left[ A^1_{\breve{\beta}} - A^2_{\breve{\beta}} \cot ( \breve{\beta} \rho) \right]}, \\[1ex]
\frac {dx_2}{dx_0} = \frac { x_2 [ A^1_{\breve{\beta}} \cot ( \breve{\beta} \rho) + A^2_{\breve{\beta}} ]}
{ \rho \left[ A^1_{\breve{\beta}} - A^2_{\breve{\beta}} \cot ( \breve{\beta} \rho) \right]}.
\end{cases}
 \end{gather*}

When $ \cos ( \breve{\beta} \rho) \neq 0$ in $\Lambda$,~\eqref{V_0-3} is represented as 
\begin{gather*}
\begin{cases}
\frac {dx_1}{dx_0} = \frac { x_1 [ A^1_{\breve{\beta}} + A^2_{\breve{\beta}} \tan ( \breve{\beta} \rho) ]}
{ \rho \left[ A^1_{\breve{\beta}} \tan ( \breve{\beta} \rho) - A^2_{\breve{\beta}} \right]}, \\[1ex]
\frac {dx_2}{dx_0} = \frac { x_2 [ A^1_{\breve{\beta}} + A^2_{\breve{\beta}} \tan ( \breve{\beta} \rho) ]}
{ \rho \left[ A^1_{\breve{\beta}} \tan ( \breve{\beta} \rho) - A^2_{\breve{\beta}} \right]}.
\end{cases}
 \end{gather*}

Relations  
\begin{gather*}
\begin{cases}
    A^1_{\breve{\beta}} \sin ( \breve{\beta} \rho) - A^2_{\breve{\beta}} \cos ( \breve{\beta} \rho) = 0,  \\[1ex]
   A^1_{\breve{\beta}} \cos ( \breve{\beta} \rho) + A^2_{\breve{\beta}} \sin ( \breve{\beta} \rho) = 0
\end{cases}
\end{gather*}
imply that zero set of $\vec V$ is empty.

$ \frac{\partial{V_{\rho}}}{\partial{x_0}} = - \breve{\beta}^2 e^{\breve{\beta} x_0}
\left[ A^1_{\breve{\beta}} \cos ( \breve{\beta} \rho) + A^2_{\breve{\beta}} \sin ( \breve{\beta} \rho) \right], \quad
\frac{\partial{V_{\rho}}}{\partial{\rho}} = \breve{\beta}^2 e^{\breve{\beta} x_0}
\left[ A^1_{\breve{\beta}} \sin ( \breve{\beta} \rho) - A^2_{\breve{\beta}} \cos ( \breve{\beta} \rho) \right]. $

The eigenvalues of~\eqref{VG tensor-merid-1} are given by formulas:
\begin{align*}
\lambda_{0}
&=  - \breve{\beta} e^{\breve{\beta} x_0}{\rho}^{-1} \left[ A^1_{\breve{\beta}} \cos ( \breve{\beta} \rho)
+ A^2_{\breve{\beta}} \sin ( \breve{\beta} \rho) \right];\\
\lambda_{1, 2}
&= \pm \breve{\beta}^2 e^{\breve{\beta} x_0} 
\sqrt{ \left[A^1_{\breve{\beta}} \cos ( \breve{\beta} \rho) + A^2_{\breve{\beta}} \sin ( \breve{\beta} \rho)\right]^2 
+ \left[ A^1_{\breve{\beta}} \sin ( \breve{\beta} \rho) - A^2_{\breve{\beta}} \cos ( \breve{\beta} \rho)\right]^2}.
\end{align*}

 Sets of degenerate points of~\eqref{VG tensor-merid-1} are provided by two independent equations:
\begin{gather*}
A^1_{\breve{\beta}} \cos ( \breve{\beta} \rho) + A^2_{\breve{\beta}} \sin ( \breve{\beta} \rho) =0;  \\[1ex]
\left[A^1_{\breve{\beta}} \cos ( \breve{\beta} \rho) + A^2_{\breve{\beta}} \sin ( \breve{\beta} \rho)\right]^2 
+ \left[ A^1_{\breve{\beta}} \sin ( \breve{\beta} \rho) - A^2_{\breve{\beta}} \cos ( \breve{\beta} \rho)\right]^2 = 0.
\end{gather*}

 \end{example}

\begin{example}

Let us  consider the radially holomorphic potentials represented by radially holomorphic fourth-order polynomials with reduced-quaternion-valued coefficients of the form 
 \begin{gather*}
 G(x) = g + I \hat{g} = \frac{I}{4} x^4 + \frac{a_1}{2} x^2, \quad a_1 = \mu_1 + I \mu_2, \quad \mu_1, \mu_2 \in \mathbb R,  \quad \mu_1 > 0.
 \end{gather*}

Generalized axially symmetric potentials $g = g(x_0,  \rho)$ are written as 
\begin{gather}
 g(x_0,  \rho) = (-x_0^3 \rho + x_0 \rho^3) + \left[ \frac{\mu_1}{2}(x_0^2 - \rho^2) - \mu_2 x_0 \rho \right].
\label{gasp-polyn-4}
\end{gather}

The radially holomorphic functions $F = F(x)$ are described as
\begin{gather*}
 F(x)= I x^3 + a_1 x.
 \end{gather*}
\begin{gather*}
 \overline{F(x)} = -I \overline{x}^3 + \overline{a_1} \overline{x}.
 \end{gather*}

The essential properties of the corresponding  axisymmetric gradient system with $1$-axial-hyperbolic harmonic velocity potential  
in $\Lambda$ $(\rho >0)$ may be demonstrated explicitly using generalized axially symmetric potentials of the form~\eqref{gasp-polyn-4}.

\begin{gather*}
V_0 = \frac{\partial{g}}{\partial{x_0}} = -3x_0^2 \rho + \rho^3 + \mu_1 x_0 - \mu_2 \rho, \qquad
V_{\rho} =  \frac{\partial{g}}{\partial{\rho}} = - x_0^3 + 3 x_0 \rho^2 - \mu_1 \rho - \mu_2 x_0.
\end{gather*}

Equation~\eqref{Petr-4} is written as  
\begin{gather*}
  \frac {dx}{dt}  
= (-3x_0^2 \rho + \rho^3 + \mu_1 x_0 - \mu_2 \rho) 
+  I (- x_0^3 + 3 x_0 \rho^2 - \mu_1 \rho - \mu_2 x_0).
  \end{gather*}

When $(-3x_0^2 \rho + \rho^3 + \mu_1 x_0 - \mu_2 \rho) \neq 0$ in $\Lambda$,~\eqref{V_0-3} is represented as 
\begin{gather*}
\begin{cases}
\frac {dx_1}{dx_0} = \frac { x_1 ( - x_0^3 + 3 x_0 \rho^2 - \mu_1 \rho - \mu_2 x_0)}{ \rho( -3x_0^2 \rho + \rho^3 + \mu_1 x_0 - \mu_2 \rho)}, \\[1ex]
\frac {dx_2}{dx_0} = \frac { x_2 ( - x_0^3 + 3 x_0 \rho^2 - \mu_1 \rho - \mu_2 x_0)}{ \rho( -3x_0^2 \rho + \rho^3 + \mu_1 x_0 - \mu_2 \rho)}.
\end{cases}
 \end{gather*}

  When $x_0 \neq 0$, relations
\begin{gather*}
\begin{cases}
V_0 =  -3x_0^2 \rho + \rho^3 + \mu_1 x_0 - \mu_2 \rho = 0, \\[1ex]
V_{\rho} = - x_0^3 + 3 x_0 \rho^2 - \mu_1 \rho - \mu_2 x_0 = 0
\end{cases}
\end{gather*}
imply that zero sets of $\vec V$ are determined from a system of equations 
\begin{gather}
\begin{cases}
  \rho^2 = x_0^2 + \mu_2, \\[1ex]
   2 x_0 \sqrt{ x_0^2 + \mu_2} = \mu_1.
\label{polyn-4}
\end{cases}
\end{gather}

Using the second equation of~\eqref{polyn-4}, we find that  
\begin{gather*}
\begin{cases}
  x_0^2 = -\frac{\mu_2}{2} \pm \frac{\sqrt{ \mu_1^2 + \mu_2^2}}{2}, \\[1ex]
  \rho^2 = \frac{\mu_2}{2} \pm \frac{\sqrt{ \mu_1^2 + \mu_2^2}}{2}.
\end{cases}
\end{gather*}

When $\mu_2 =0$, zero sets of $\vec V$ are determined from a system of equations 
\begin{gather*}
\begin{cases}
  \rho^2 = x_0^2, \\[1ex]
  2 x_0^2 = \mu_1.
\end{cases}
\end{gather*}

\begin{gather*}
\frac{\partial{V_{\rho}}}{\partial{x_0}} = - 3 x_0^2 + 3 \rho^2 - \mu_2, \qquad
\frac{\partial{V_{\rho}}}{\partial{\rho}} = 6 x_0 \rho -\mu_1.
\end{gather*}

The eigenvalues of~\eqref{VG tensor-merid-1} are given by formulas:
\begin{gather*}
\lambda_{0} = - \frac{x_0^3}{\rho} + 3 x_0 \rho - \mu_1 - \mu_2 \frac{x_0}{\rho}; \\[1ex]
\lambda_{1, 2} = \pm \sqrt{ ( 3 x_0^2 - 3 \rho^2 + \mu_2)^2 + (6 x_0 \rho -\mu_1)^2}.
\end{gather*}

 Sets of degenerate points of~\eqref{VG tensor-merid-1} are provided by two independent equations:
\begin{align*}
  x_0^3 - 3 x_0 \rho^2 + \mu_1 \rho + \mu_2 x_0 =0; \\
 ( 3 x_0^2 - 3 \rho^2 + \mu_2)^2 + (6 x_0 \rho -\mu_1)^2 = 0.
\end{align*}

\end{example}

\section{Concluding Remarks}
 \label{sec6}
 Consider basic properties of smooth gradient systems with scalar potential $h$ in the phase space $\mathbb R^4=\{(x_0, x_1,x_2, x_3)\}$:
 \begin{gather}
\begin{cases}
\frac {dx_0}{dt} = V_0(x_0,x_1,x_2,x_3) = \frac{\partial{h(x_0,x_1,x_2,x_3)}}{\partial{x_0}}, \\[1ex]
\frac {dx_1}{dt} = V_1(x_0,x_1,x_2,x_3) = \frac{\partial{h(x_0,x_1,x_2,x_3)}}{\partial{x_1}}, \\[1ex]
\frac {dx_2}{dt} = V_2(x_0,x_1,x_2,x_3) = \frac{\partial{h(x_0,x_1,x_2,x_3)}}{\partial{x_2}}, \\[1ex]
\frac {dx_3}{dt} = V_3(x_0,x_1,x_2,x_3) = \frac{\partial{h(x_0,x_1,x_2,x_3)}}{\partial{x_3}}.
\end{cases}
    \label{traject-4}
 \end{gather}
The vector $\vec V = (V_0, V_1, V_2, V_3)$ implies a rich variety of analytic models of four-dimensional potential  vector fields. 
Zero sets of $\vec V$ in simply connected open domains $\Lambda \subset \mathbb R^4$
coincide with the critical sets of the scalar potential $h$ in $\Lambda$.

Consider the following four-dimensional extension of the system~\eqref{alpha-axial-hyperbolic-system-3}:
\begin{gather}
\begin{cases}
  (x_1^2+x_2^2+x_3^2) \ \mathrm{div}\ { \vec V} -
 \alpha \left(x_1 V_1 + x_2 V_2 + x_3 V_3 \right) =0,  \\[1ex]
    \frac{\partial{V_0}}{\partial{x_1}}= \frac{\partial{V_1}}{\partial{x_0}},
       \quad \frac{\partial{V_0}}{\partial{x_2}}= \frac{\partial{V_2}}{\partial{x_0}},
       \quad \frac{\partial{V_0}}{\partial{x_3}}= \frac{\partial{V_3}}{\partial{x_0}}, \\[1ex]
      \frac{\partial{V_1}}{\partial{x_2}}= \frac{\partial{V_2}}{\partial{x_1}},
       \quad \frac{\partial{V_1}}{\partial{x_3}}= \frac{\partial{V_3}}{\partial{x_1}},
       \quad \frac{\partial{V_2}}{\partial{x_3}}= \frac{\partial{V_3}}{\partial{x_2}}.
\end{cases}
\label{alpha-axial-isotropic-system-4}
\end{gather}

The scalar potential $h$ in $\Lambda$,  where $\vec V = \mathrm{grad} \ h,$ 
allows us to reduce every $C^1$-solution 
of the system~\eqref{alpha-axial-isotropic-system-4} to a $C^2$-solution of a second-order elliptic equation with three singular coefficients:
 \begin{gather}
(x_1^2+ x_2^2 + x_3^2)\Delta{h} - \alpha \left( x_1\frac{\partial{h}}{\partial{x_1}} +
x_2\frac{\partial{h}}{\partial{x_2}} + x_3\frac{\partial{h}}{\partial{x_3}} \right)  =0.
 \label{alpha-axial-hyperbolic-4}
  \end{gather}

The four-dimensional potential vector field $\vec V = V(x_0, x_1,x_2,x_3)$ 
in the phase space $\mathbb R^4=\{(x_0, x_1,x_2, x_3)\}$ may be seen as the gradient vector field, 
where $\vec V(\vec x) = -\mathrm{grad} \ \hslash(\vec x)$, $\hslash = \hslash(\vec x) $ denotes a potential function. 
The coefficient $\phi = \phi(x_0,x_1,x_2,x_3)$ may be seen as the density of the integral invariant
 $\int_{\Lambda} \phi(x_0,x_1,x_2,x_3)dx_0 dx_1 dx_2 dx_3$, respectively 
(see, e.g., \cite{Poincare:1892-1993,NemStep:1960,Goriely:2001,Perko:2001,Wiggins:2003,Strogatz:2018}).

The first-order system 
\begin{gather}
\begin{cases}
 (x_1^2+x_2^2+x_3^2) \left( \frac{\partial{u_0}}{\partial{x_0}}-
      \frac{\partial{u_1}}{\partial{x_1}}-\frac{\partial{u_2}}{\partial{x_2}}-
      \frac{\partial{u_3}}{\partial{x_3}} \right) + \alpha (x_1u_1+x_2u_2+x_3u_3)=0 \\[1ex]
      \frac{\partial{u_0}}{\partial{x_1}}=-\frac{\partial{u_1}}{\partial{x_0}},
       \quad \frac{\partial{u_0}}{\partial{x_2}}=-\frac{\partial{u_2}}{\partial{x_0}},
       \quad \frac{\partial{u_0}}{\partial{x_3}}=-\frac{\partial{u_3}}{\partial{x_0}}, \\[1ex]
      \frac{\partial{u_1}}{\partial{x_2}}=\ \ \frac{\partial{u_2}}{\partial{x_1}},
       \quad \frac{\partial{u_1}}{\partial{x_3}}=\ \ \frac{\partial{u_3}}{\partial{x_1}},
       \quad \frac{\partial{u_2}}{\partial{x_3}}=\ \ \frac{\partial{u_3}}{\partial{x_2}},
\end{cases}
\label{eq:A_4^alpha-system}
\end{gather}
where  $(u_0, u_1, u_2, u_3) :=(V_0, -V_1, -V_2, -V_3)$, demonstrates explicitly 
 a family of axially symmetric generalizations of the Cauchy-Riemann system  in $\mathbb R^4$
in accordance with the system~\eqref{alpha-axial-isotropic-system-4} for different values of the parameter $\alpha$
(see \cite{Br:Hefei2020,Br:2021}).

Analytic models of four-dimensional harmonic potential fields $\vec V$ 
satisfy the Riesz system in $\mathbb R^4$  (see, e.g., \cite{SteinWeiss:1971})
\begin{gather*}
\begin{cases}
       \frac{\partial{V_0}}{\partial{x_0}} + \frac{\partial{V_1}}{\partial{x_1}} + \frac{\partial{V_2}}{\partial{x_2}} + \frac{\partial{V_3}}{\partial{x_3}} =0,  \\[1ex]
    \frac{\partial{V_0}}{\partial{x_1}}= \frac{\partial{V_1}}{\partial{x_0}},
       \quad \frac{\partial{V_0}}{\partial{x_2}}= \frac{\partial{V_2}}{\partial{x_0}},
       \quad \frac{\partial{V_0}}{\partial{x_3}}= \frac{\partial{V_3}}{\partial{x_0}}, \\[1ex]
      \frac{\partial{V_1}}{\partial{x_2}}= \frac{\partial{V_2}}{\partial{x_1}},
       \quad \frac{\partial{V_1}}{\partial{x_3}}= \frac{\partial{V_3}}{\partial{x_1}},
       \quad \frac{\partial{V_2}}{\partial{x_3}}=
       \frac{\partial{V_3}}{\partial{x_2}}.
\end{cases}
\end{gather*}
General class of exact solutions of the Riesz system in $\mathbb R^4$ 
in the context of \emph{Quaternionic analysis in $\mathbb R^4$}
(see, e.g., \cite{BraDel:2003,LeZe:CMFT2004,Del:2007})
is equivalently represented as general class of analytic solutions of the system
\begin{gather*}
(R)
\begin{cases}
 \frac{\partial{u_0}}{\partial{x_0}}-
      \frac{\partial{u_1}}{\partial{x_1}}- \frac{\partial{u_2}}{\partial{x_2}} - \frac{\partial{u_3}}{\partial{x_3}} =0,  \\[1ex]
       \frac{\partial{u_0}}{\partial{x_1}}=-\frac{\partial{u_1}}{\partial{x_0}},
     \quad \frac{\partial{u_0}}{\partial{x_2}}=-\frac{\partial{u_2}}{\partial{x_0}},
       \quad \frac{\partial{u_0}}{\partial{x_3}}=-\frac{\partial{u_3}}{\partial{x_0}}, \\[1ex]
      \frac{\partial{u_1}}{\partial{x_2}}=\ \ \frac{\partial{u_2}}{\partial{x_1}},
       \quad \frac{\partial{u_1}}{\partial{x_3}}=\ \ \frac{\partial{u_3}}{\partial{x_1}},
       \quad \frac{\partial{u_2}}{\partial{x_3}}=\ \ \frac{\partial{u_3}}{\partial{x_2}}.
\end{cases}
\end{gather*}
Exact solutions of the system $(R)$ in $\mathbb R^4$ are referred to as the quaternion-valued monogenic functions 
$u= u_0 + iu_1 + ju_2 + ku_3$  with harmonic components $u_l= u_l(x_0,x_1,x_2,x_3)$ $(l= 0,1,2,3)$.

When $\alpha > 0$, the system~\eqref{eq:A_4^alpha-system} may be characterized as 
$\alpha$-axial-hyperbolic non-Euclidean modification of the system $(R)$ with respect to the conformal
metric defined outside the axis $x_0$ by the formula
$$
ds^2 = \frac{d{x_0}^2 + d{x_1}^2 + d{x_2}^2 +
d{x_3}^2}{\rho^{\alpha}}.
$$

\begin{definition}
Every exact solution of eqn~\eqref{alpha-axial-hyperbolic-4} under the condition $\alpha>0$ 
 in a simply connected open domain $\Lambda \subset \mathbb R^4$ $(\rho > 0)$
 is called $\alpha$-axial-hyperbolic harmonic potential in $\Lambda$.
\end{definition}

The analytic properties of exact solutions of system~\eqref{eq:A_4^alpha-system} may be studied in more detail 
in the context of \emph{Quaternionic analysis in $\mathbb R^4$ and its non-Euclidean modifications}
(see, e.g., \cite{HempLeut:1996,LeZe:CMFT2004,ErOrel:2019,Br:Hefei2020}) using the well-known formulas \cite{Sudbery:1979}
  \begin{gather*}
   \begin{cases}
  x_0 = \frac{1}{4}(x - ixi - jxj - kxk), \\[1ex]
   x_1 = \frac{1}{4i}(x - ixi + jxj + kxk), \\[1ex]
   x_2 = \frac{1}{4j}(x + ixi - jxj + kxk), \\[1ex]
   x_3 = \frac{1}{4k}(x + ixi + jxj - kxk).
   \end{cases}
   \end{gather*}
   
\bmhead{Acknowledgments} 
The author would like to thank Prof. Aksenov~A.V., Prof. Astashova~I.V., Prof. Filinovskiy~A.V., Prof. K\"{a}hler~U., Prof. Khesin~B.A., 
Dr. Kisil~V.V., Prof. Kravchenko~V.V., Prof. Leutwiler~H., Prof. Sergeev~I.N., Prof. Spr\"{o}{\ss}ig~W., Prof. Urinov~A.K.
 for inspiring discussions, that have led to the emergence of this paper.

\end{document}